# NONPARAMETRIC REGRESSION PENALIZING DEVIATIONS FROM ADDITIVITY[1]

By M. Studer, B. Seifert and T. Gasser

*University of Zurich*

Due to the curse of dimensionality, estimation in a multidimensional nonparametric regression model is in general not feasible. Hence, additional restrictions are introduced, and the additive model takes a prominent place. The restrictions imposed can lead to serious bias. Here, a new estimator is proposed which allows penalizing the non-additive part of a regression function. This offers a smooth choice between the full and the additive model. As a byproduct, this penalty leads to a regularization in sparse regions. If the additive model does not hold, a small penalty introduces an additional bias compared to the full model which is compensated by the reduced bias due to using smaller bandwidths.

For increasing penalties, this estimator converges to the additive smooth backfitting estimator of Mammen, Linton and Nielsen [*Ann. Statist.* **27** (1999) 1443–1490].

The structure of the estimator is investigated and two algorithms are provided. A proposal for selection of tuning parameters is made and the respective properties are studied. Finally, a finite sample evaluation is performed for simulated and ozone data.

**1. Introduction.** Let $(\underline{\mathbf{X}}_i, \varepsilon_i)$, $i = 1, \ldots, n$, be independent identically distributed random vectors with $\underline{\mathbf{X}}_i \in [0,1]^d$. Define the response as $Y_i = r^{\text{true}}(\underline{\mathbf{X}}_i) + \varepsilon_i$. The errors $\varepsilon_i$ have expectation zero and variance $\sigma^2$ and are independent of $\underline{\mathbf{X}}_i$. The goal is to estimate $r^{\text{true}}(\underline{\mathbf{x}})$ given data $(\underline{\mathbf{X}}_i, Y_i)$.

In the *full model*, we assume only that the unknown regression function

$$r^{\text{true}}(\underline{\mathbf{x}}) = \mathbb{E}(Y | \underline{\mathbf{X}} = \underline{\mathbf{x}})$$

is smooth. Specifically, we assume that $r^{\text{true}}$ is twice continuously differentiable as we will use a local linear estimator. The rate of convergence of mean square error is $O(n^{-4/(4+d)})$ [Stone (1980, 1982)].

Received June 2003; revised June 2004.

[1]Supported by Swiss National Science Foundation Project 2100-052567.97.

Estimating in the full model suffers from the "curse of dimensionality." This leads to consideration of less general models. In the *additive model* it is assumed that the regression function has the special form

$$(1) \qquad r^{\text{true}}(\underline{\mathbf{x}}) = r^{\text{true}}_{\text{add},0} + r^{\text{true}}_{\text{add},1}(x_1) + \cdots + r^{\text{true}}_{\text{add},d}(x_d).$$

The rate of convergence is $O(n^{-4/5})$ as for $d = 1$ [Stone (1985, 1986)].

Choosing the additive model may lead to serious bias due to neglecting the nonadditive component of the regression function. Estimating the full model may, however, lead also to a large bias since a large (optimal) bandwidth has to be used to achieve the same rate for variance as for squared bias.

In this paper we introduce a parametric family of estimators $\widehat{\underline{r}}_R$ ($R \geq 0$) which includes asymptotically optimal estimators for the full ($R = 0$) and the additive ($R = \infty$) model as *special cases*. The aim is to offer a *continuous model choice* via the tuning parameter $R$.

The philosophy behind additive modeling might be described as follows: rather than assuming the strict validity of the additive assumption, one goes for the additive part of the underlying regression function to avoid the curse of dimensionality. The approach of this paper offers us more flexibility in case of highly nonadditive functions: instead of switching to the full model (or tolerating a large bias for the additive fit), one chooses a fit in between, which takes into account part of the nonadditive structure.

*Local linear estimation.* For fixed $\underline{\mathbf{x}}$, let $\widehat{\underline{\beta}} = (\widehat{\beta}_0, \ldots, \widehat{\beta}_d)$ be the minimizer of

$$(2) \qquad SSR(\underline{\beta}, \underline{\mathbf{x}}) = \frac{1}{n} \sum_{i=1}^{n} \left( Y_i - \beta_0 - \sum_{k=1}^{d} \beta_k \frac{X_{i,k} - x_k}{h_k} \right)^2 K_h(\underline{\mathbf{X}}_i, \underline{\mathbf{x}}),$$

where $K_h(\underline{\mathbf{X}}_i, \underline{\mathbf{x}}) = K(\text{diag}(h_1, \ldots, h_d)^{-1}(\underline{\mathbf{X}}_i - \underline{\mathbf{x}}))/(h_1 \cdots \cdots h_d) \geq 0$ is the kernel weight of the observation $(\underline{\mathbf{X}}_i, Y_i)$ for the output point $\underline{\mathbf{x}}$. The bandwidths $h_1, \ldots, h_d$ are scale parameters. We assume that $h_1, \ldots, h_d$ are of the same order and set $h = \sqrt[d]{h_1 \cdots \cdots h_d}$. The diagonal matrix with diagonal elements $h_1, \ldots, h_d$ is denoted by $\text{diag}(h_1, \ldots, h_d)$. The local linear estimator of $r^{\text{true}}(\underline{\mathbf{x}})$ at output point $\underline{\mathbf{x}}$ is $\widehat{\beta}_0$.

Under usual regularity conditions, variance is proportional to $(nh^d)^{-1}$ and squared bias is proportional to $h^4$. The optimal rate for the MSE is $n^{-4/(4+d)}$, using a bandwidth $h$ proportional to $n^{-1/(4+d)}$. The local linear estimator achieves asymptotically the linear minimax risk when using spherically symmetric Epanechnikov kernels. This optimality result was shown in Fan (1993) for $d = 1$ and in Fan, Gasser, Gijbels, Brockmann and Engel (1997) for $d > 1$. For finite sample size, however, regularization is an issue [Seifert and Gasser (1996)], as the variance is unbounded in sparse regions.

As we will see later, our modeling approach via the parameter $R$ leads, as a byproduct, also to a regularization.

Mammen, Linton and Nielsen (1999) (referred to as MLN below) introduced a backfitting estimator for the additive model which achieves the same asymptotics as the *oracle estimator*, which is a univariate local linear estimator for data $(X_{i,k}, Y_i - \sum_{\kappa \neq k} r_{\mathrm{add},\kappa}^{\mathrm{true}}(X_{i,\kappa}))$. Consequently, it inherits the above mentioned optimality. They evaluated a local linear estimator on a continuum (e.g., $[0,1]^d$), using a vector of parameter functions $\underline{r}(\underline{\mathbf{x}}) = (r^0(\underline{\mathbf{x}}), \ldots, r^d(\underline{\mathbf{x}}))$. The first function $r^0(\underline{\mathbf{x}})$ is the *intercept* (i.e., $\beta_0$ for $\underline{\mathbf{x}}$) and the other functions are *slopes* (i.e., $\beta_1, \ldots, \beta_d$). MLN decompose $\underline{r}(\underline{\mathbf{x}})$ into additive ($\underline{r}_{\mathrm{add}}$) and orthogonal ($\underline{r}_\perp$) components, and set the orthogonal component to zero:

$$\widehat{\underline{r}}_{\mathrm{add}} = \arg\min_{\underline{r}_{\mathrm{add}}} \int_{[0,1]^d} SSR(\underline{r}_{\mathrm{add}}(\underline{\mathbf{x}}), \underline{\mathbf{x}}) \, d\underline{\mathbf{x}}.$$

The estimator $\widehat{\underline{r}}_{\mathrm{add}}$ has an interpretation as a projection $\mathcal{P}_* \widehat{\underline{r}}_{ll}$ of the local linear $\widehat{\underline{r}}_{ll}$ to the additive subspace.

Instead of a projection we use a *penalty* $R$ to shrink the orthogonal component towards zero. Formally,

$$\widehat{\underline{r}}_R = \arg\min_{\underline{r}} \int_{[0,1]^d} SSR(\underline{r}(\underline{\mathbf{x}}), \underline{\mathbf{x}}) \, d\underline{\mathbf{x}} + R \|\underline{r}_\perp\|_2^2.$$

For $R = 0$ we get the usual local linear estimator, and for $R = \infty$ we obtain the additive estimator of MLN. For general $R$ we get a family of estimators connecting $\widehat{\underline{r}}_{ll}$ with $\widehat{\underline{r}}_{\mathrm{add}}$ with common additive part $\mathcal{P}_* \widehat{\underline{r}}_R = \widehat{\underline{r}}_{\mathrm{add}}$.

EXAMPLE. Let us now illustrate the benefit of a smooth choice between full and additive models for some simulated data with known regression function and random uniform design; see Figure 1.

Originally this realization of the data was used in Seifert and Gasser (2000) in the context of locally ridging the local linear estimator. (Another 50 realizations are summarized in Section 5.1.) Due to symmetry of the true regression function $[r^{\mathrm{true}}(x_1, x_2) = r^{\mathrm{true}}(x_2, x_1)]$, there is no need to consider separate bandwidths for each coordinate. Note that the smoothing windows have the same size in the interior and at the boundary by choosing a larger bandwidth at the boundary [see Figure 1(a)]. We use a product Epanechnikov kernel. The output grid consists of $50 \times 50$ points and the parameters are the minimizers of integrated squared residuals (ISE); see Figure 3(a).

Even though the regression function is clearly nonadditive, penalizing the nonadditive part leads to a remarkable improvement in optimal ISE from 8.3 [Figure 2(a)] to 6.0 [Figure 2(d)]. A small penalty $R$ stabilizes output points where the local linear estimator is wiggly but has little effect in well-determined regions [Figure 2(b) vs. 2(d)]. This illustrates another useful property of penalizing: regularization of the local linear estimator.

Generalizing a method often improves goodness of fit, while parameter selection becomes more difficult. Let us apply $AIC_C$ for selecting parameters $R$ and $h$ [Hurvich, Simonoff and Tsai (1998)]. This criterion tries to find a compromise between good fit and small complexity of the model (i.e., low trace of hat matrix). Figure 3 shows that parameter selection is successful in this example.

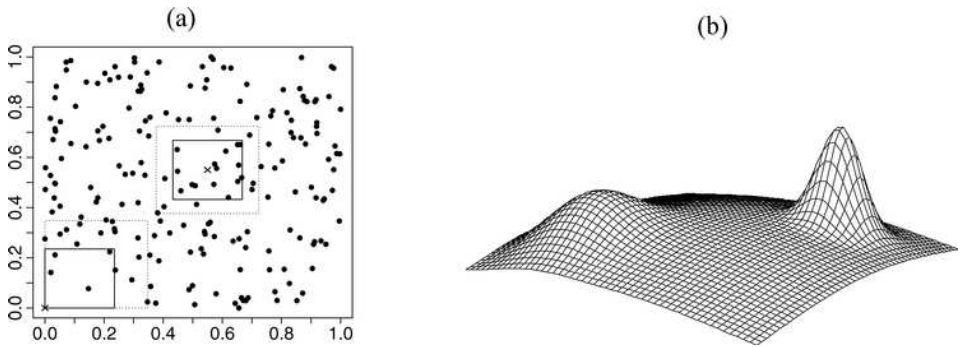

FIG. 1.   *Simulated data using $n = 200$ random observations* (a) *(design $n = 200$) and regression function* (b) *(true regression function) (range $[9, 54]$, residual variance $\sigma^2 = 25$). Smoothing windows are of constant size due to increased bandwidth at the boundary:* (a) *displays smoothing windows for $h = 0.117$ and $h = 0.174$ at output points $(0.55, 0.55)$ and $(0, 0)$.*

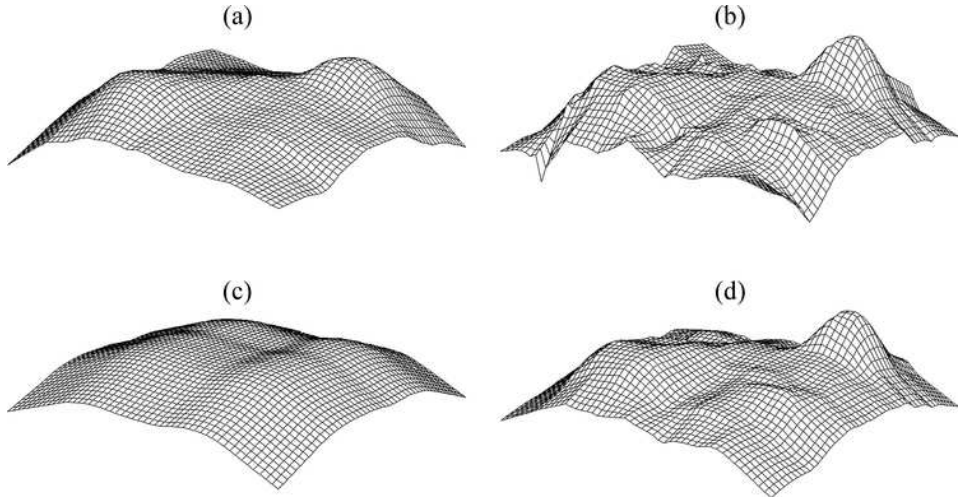

FIG. 2.   *Comparison of different estimators. The local linear estimator is either heavily biased* (a) *($h = 0.174$, $R = 0$, ISE $= 8.3$) or wiggly* (b) *($h = 0.117$, $R = 0$, ISE $= 8.8$). Additive estimation* (c) *($h = 0.197$, $R = \infty$, ISE $= 17$) is even worse. The penalized estimator* (d) *($h = 0.117$, $R = 0.163$, ISE $= 6.0$) is stabilized without oversmoothing: ISE is improved by more than a quarter.*

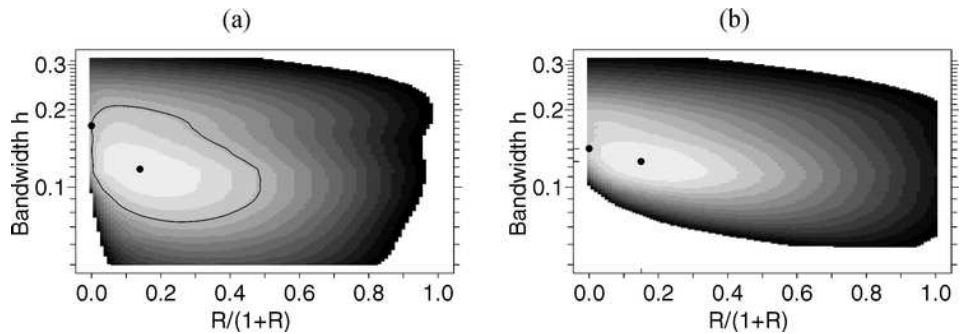

Fig. 3. *Comparison of ISE* (a) [*integrated squared error (ISE)*] *and* AIC_C (b) *as a function of bandwidth h (log-scale) and penalty R ($\frac{R}{1+R}$-scale). The global minimum of ISE is at (h = 0.117, R = 0.163). A contour line bounds a region of parameters outperforming the ordinary local linear estimator (minimum at h = 0.174).*

*Contents.* The paper is organized as follows: Section 2 defines the proposed estimator both in a discrete and in a continuous version. A computationally efficient direct and an iterative algorithm are developed in Section 3. Properties of the estimator are studied in Section 4: the penalized estimator is shown to be a pointwise compromise between an additive and the local linear fit. A decomposition into an additive part and an orthogonal remainder term is derived, where only the nonadditive part involves shrinking. In addition to model flexibility, the approach offers a regularization in sparse regions of the design. We then justify the interpretation of the local linear and the additive estimators as special cases of $\widehat{\underline{r}}_R$ for $R \to 0$ or $\infty$. The convergence of $\widehat{\underline{r}}_R$ to the MLN estimator for $R \to \infty$ is investigated. The data-adaptive simultaneous choice of $h$ and $R$ is analyzed to some extent. Section 5 is devoted to a simulation study. Furthermore, the estimator is applied to the ozone dataset. A summary of the contents is provided at the beginning of each section. Software is available on our homepage www.biostat.unizh.ch.

**2. Definition of the penalized estimator.** A local linear estimator is evaluated on a set of output points. For penalizing deviations from the additive model, these output points should form a product space. In Section 2.1 we choose the interval $[0, 1]^d$ as a continuous set of output points and start with definitions similar to MLN. This choice is suitable for deriving theoretical properties. In practice the continuous set of output points is approximated by an equidistant grid,

$$\left\{0, \frac{1}{m_1-1}, \frac{2}{m_1-1}, \ldots, 1\right\} \times \cdots \times \left\{0, \frac{1}{m_d-1}, \frac{2}{m_d-1}, \ldots, 1\right\}$$

as in Section 2.2.

2.1. *Estimation on an interval.* We will introduce a Hilbert space $(\mathcal{F}, \|\cdot\|_*)$ such that the local linear estimator $\widehat{\underline{r}}_{ll}$ corresponds to a projection of the response $\underline{Y}$ to some subspace $\mathcal{F}_{\text{full}} \subset \mathcal{F}$.

MLN consider a subspace $\mathcal{F}_{\text{add}} \subset \mathcal{F}_{\text{full}}$ of additive functions and obtain $\widehat{\underline{r}}_{\text{add}}$ by projecting $\underline{Y}$ to $\mathcal{F}_{\text{add}}$.

We consider another norm $\|\cdot\|_R$ being the sum of $\|\cdot\|_*$ and some penalty with parameter $R$ on the squared distance from $\mathcal{F}_{\text{add}}$. The penalized estimator $\widehat{\underline{r}}_R$ is the projection with respect to $\|\cdot\|_R$ of $\underline{Y}$ to $\mathcal{F}_{\text{full}}$.

Define the vector space of $(n+1)(d+1)$ functions

$$\mathcal{F} = \{\underline{r} = (r^{i,\ell}|i=0,\ldots,n; \ell=0,\ldots,d)|r^{i,\ell} \colon [0,1]^d \to \mathbb{R}\}.$$

Let us define the projection $\mathcal{P}_0$ on $\mathcal{F}$, which replaces $r^{i,\ell}$ by $r^{0,\ell}$. In other words, if $\underline{\breve{r}} = \mathcal{P}_0\underline{r}$, then $\breve{r}^{i,\ell}(\underline{\mathbf{x}}) = r^{0,\ell}(\underline{\mathbf{x}})$. The image of $\mathcal{P}_0$ is denoted by $\mathcal{F}_{\text{full}}$. For simplicity of notation, the index $i$ is omitted:

$$\mathcal{F}_{\text{full}} = \{\underline{r} = (r^0, \ldots, r^d)|r^\ell \colon [0,1]^d \to \mathbb{R}, \ell = 0, \ldots, d\}.$$

The observations $Y_i$, $i = 1, \ldots, n$, are coded as $\underline{r}_Y \in \mathcal{F}$ by

$$r_Y^{i,\ell}(\underline{\mathbf{x}}) = \begin{cases} Y_i, & \text{for } i > 0 \text{ and } \ell = 0, \\ 0, & \text{otherwise.} \end{cases}$$

Define the design-dependent seminorm $\|\cdot\|_*$ on $\mathcal{F}$ by

$$\|\underline{r}\|_*^2 = \int \frac{1}{n}\sum_{i=1}^n \left[r^{i,0}(\underline{\mathbf{x}}) + \sum_{k=1}^d r^{i,k}(\underline{\mathbf{x}})\frac{X_{i,k}-x_k}{h_k}\right]^2 K_h(\underline{\mathbf{X}}_i, \underline{\mathbf{x}})\, d\underline{\mathbf{x}},$$

where $K_h(\underline{\mathbf{X}}_i, \underline{\mathbf{x}})$ is the kernel weight of the observation $(\underline{\mathbf{X}}_i, Y_i)$ for the output point $\underline{\mathbf{x}}$.

Hence, for $\underline{r} \in \mathcal{F}_{\text{full}}$ we have

$$(3) \quad \|\underline{r}_Y - \underline{r}\|_*^2 = \int \frac{1}{n}\sum_{i=1}^n \left[Y_i - r^0(\underline{\mathbf{x}}) - \sum_{k=1}^d r^k(\underline{\mathbf{x}})\frac{X_{i,k}-x_k}{h_k}\right]^2 K_h(\underline{\mathbf{X}}_i, \underline{\mathbf{x}})\, d\underline{\mathbf{x}}$$

and the integrand corresponds to the minimization problem for the local linear estimator. Consequently, we denote the minimizer by $\widehat{\underline{r}}_{ll}$.

The interpretation of $\widehat{\underline{r}}_{ll}$ as projection of $\underline{r}_Y$ to $\mathcal{F}_{\text{full}}$ was developed by Mammen, Marron, Turlach and Wand (2001) and is quite useful when incorporating constraints, that is, minimizing (3) for $\underline{r}$ in a subset of $\mathcal{F}_{\text{full}}$.

Consider now an additive subspace $\mathcal{F}_{\text{add}} \subset \mathcal{F}_{\text{full}}$:

$$\mathcal{F}_{\text{add}} = \{\, \underline{r} \in \mathcal{F}_{\text{full}} | r^0(\underline{\mathbf{x}}) \text{ is additive};$$

$$\text{for } k = 1, \ldots, d,\ r^k(\underline{\mathbf{x}}) \text{ depends only on } x_k \}.$$

Define the additive estimator $\widehat{\underline{r}}_{\text{add}}$ as the minimizer for $\underline{r} \in \mathcal{F}_{\text{add}}$ of $\|\underline{r}_Y - \underline{r}\|_*^2$.

Projecting a Nadaraya–Watson estimator to an additive subspace was first considered by Nielsen and Linton (1998) for $d = 2$ and extended to higher dimensions in MLN. The projected local linear estimator $\widehat{\underline{r}}_{\text{add}}$ was introduced by MLN and has attractive properties. Nielsen and Sperlich (2005) discuss practical aspects of this estimator, which is called smooth backfitting there. These include implementation, parameter selection by cross validation and finite sample evaluation.

Let us introduce further notation. Define the $L^2$-norm on $\mathcal{F}$ by

$$\|\underline{r}\|_2^2 = \frac{1}{n+1} \sum_{i=0}^{n} \sum_{\ell=0}^{d} \int [r^{i,\ell}(\underline{\mathbf{x}})]^2 \, d\underline{\mathbf{x}}.$$

Denote by $\mathcal{P}_{\text{add}}$ the $\|\cdot\|_2$-orthogonal projection from $\mathcal{F}_{\text{full}}$ into $\mathcal{F}_{\text{add}}$.

More formally, we define $\underline{r}_{\text{add}} = \mathcal{P}_{\text{add}} \underline{r}$ via $r_{\text{add}}^0(\underline{\mathbf{x}}) = \sum_{k=1}^{d} \int r^0(\underline{\mathbf{x}}) \, d\underline{\mathbf{x}}_{-k} - (d-1) \int r^0(\underline{\mathbf{x}}) \, d\underline{\mathbf{x}}$ and $r_{\text{add}}^k(\underline{\mathbf{x}}) = \int r^k(\underline{\mathbf{x}}) \, d\underline{\mathbf{x}}_{-k}$, where $\int \cdots d\underline{\mathbf{x}}_{-k}$ denotes the integral with respect to all components of $\underline{\mathbf{x}}$ except $x_k$. Furthermore, let $\mathcal{P}_*$ be the $\|\cdot\|_*$ projection from $\mathcal{F}$ (or $\mathcal{F}_{\text{full}}$) to $\mathcal{F}_{\text{add}}$; see Appendix A.0.10.

Next, a penalty on the nonadditive part of $\underline{r}$ is added to $\|\cdot\|_*$. Define the seminorm $\|\cdot\|_R$ on $\mathcal{F}$:

$$\|\underline{r}\|_R^2 = \|\underline{r}\|_*^2 + R\|(\mathcal{I} - \mathcal{P}_{\text{add}})\mathcal{P}_0 \underline{r}\|_2^2,$$

where $\mathcal{I}$ is the identity. The penalized estimator $\widehat{\underline{r}}_R$ is defined as the minimizer of

$$\|\underline{r}_Y - \underline{r}\|_R^2 = \int \frac{1}{n} \sum_{i=1}^{n} \left[ Y_i - r^0(\underline{\mathbf{x}}) - \sum_{k=1}^{d} r^k(\underline{\mathbf{x}}) \frac{X_{i,k} - x_k}{h_k} \right]^2 K_h(\underline{\mathbf{X}}_i, \underline{\mathbf{x}}) \, d\underline{\mathbf{x}}$$

$$(4) \qquad\qquad + R\|(\mathcal{I} - \mathcal{P}_{\text{add}})\underline{r}\|_2^2$$

for $\underline{r} \in \mathcal{F}_{\text{full}}$. For the penalty term we use the fact that $\mathcal{P}_0$ is the identity on $\mathcal{F}_{\text{full}}$ and that $\mathcal{P}_0 \underline{r}_Y = 0$. The latter was the reason for introducing the components with $i = 0$ in $\mathcal{F}$. Properties of $\widehat{\underline{r}}_R$ will be analyzed in Section 4.

REMARK ON THE CHOICE OF THE PENALTY IN (4). For any choice of the penalty, the MLN estimator $\widehat{\underline{r}}_{\text{add}}$ would be the additive part of $\widehat{\underline{r}}_R$ with respect to the norm $\|\cdot\|_*$, assuming invariance under addition of an element $\mathcal{F}_{\text{add}}$ to $\underline{r}$ in (4); see Proposition 4 in Section 4.2.

We used $\|\cdot\|_2$ for the penalty instead of $\|\cdot\|_*$ because the latter is inferior in *sparse regions*. Moreover, $\mathcal{P}_*$ is not self-adjoint with respect to $\|\cdot\|_2$.

2.2. *Estimation on a grid.* Now an approximation of (4) on a grid is derived.

Let the output grid

$$\{t_1^1, \ldots, t_{m_1}^1\} \times \{t_1^2, \ldots, t_{m_2}^2\} \times \cdots \times \{t_1^d, \ldots, t_{m_d}^d\} \subset [0,1]^d$$

consist of $m_k$ values for the $k$th coordinate and enumerate its $m = m_1 \times \cdots \times m_d$ output points by $\underline{\mathbf{t}}_j = (t_{j,1}, \ldots, t_{j,d})$ for $j = 1, \ldots, m$. In order to get an appropriate approximation of (4), the output grid has to be sufficiently dense and has to increase with $n$. Denote by $\underline{\beta}^{(j)} \in \mathbb{R}^{d+1}$ the parameters of the local linear estimator at $\underline{\mathbf{t}}_j$. The parameter space for the local linear estimator on the output grid is

$$\mathbf{F}_{\text{full}} = \{\underline{\boldsymbol{\beta}} = \text{col}_j(\underline{\boldsymbol{\beta}}^{(j)}) | \underline{\boldsymbol{\beta}}^{(j)} \in \mathbb{R}^{d+1}\} = \mathbb{R}^{m(d+1)},$$

where $\text{col}_j(\underline{\boldsymbol{\beta}}^{(j)})$ denotes the column vector obtained by vertically stacking $\underline{\boldsymbol{\beta}}^{(1)}, \ldots, \underline{\boldsymbol{\beta}}^{(m)}$. The accompanying norm is the Euclidean norm $\|\cdot\|$. The additive subspace is defined as

$$\mathbf{F}_{\text{add}} = \{\text{col}_j(\underline{\boldsymbol{\beta}}^{(j)}) | \exists \underline{r}_{\text{add}} \in \mathcal{F}_{\text{add}} : \beta_\ell^{(j)} = r_{\text{add}}^\ell(\underline{\mathbf{t}}_j)\}.$$

Let $\mathbf{P}_{\text{add}}$ be the orthogonal projection from $\mathbf{F}_{\text{full}}$ to $\mathbf{F}_{\text{add}}$. The local linear estimator $\underline{\widehat{\boldsymbol{\beta}}}_{ll}^{(j)}$ at output point $\underline{\mathbf{t}}_j$ is the minimizer of the sum of weighted squared residuals SSR [see (2) in Section 1]. The simultaneous local linear estimator on the grid minimizes the sum of SSR over all output points $\underline{\mathbf{t}}_j$. Finally, we add a penalty proportional to the squared distance of the parameters to the additive submodel,

$$(5) \qquad \underline{\widehat{\boldsymbol{\beta}}}_R = \underset{\underline{\boldsymbol{\beta}} \in \mathbf{F}_{\text{full}}}{\arg \min} \sum_{j=1}^m SSR(\underline{\boldsymbol{\beta}}^{(j)}, \underline{\mathbf{t}}_j) + R\|(\mathbf{I} - \mathbf{P}_{\text{add}})\underline{\boldsymbol{\beta}}\|^2.$$

The penalized estimator $\widehat{r}_R$ is the intercept of $\underline{\widehat{\boldsymbol{\beta}}}_R^{(j)}$, that is, $\widehat{r}_R(\underline{\mathbf{t}}_j) = [\underline{\widehat{\boldsymbol{\beta}}}_R^{(j)}]_0$.

An efficient algorithm will be presented in Section 3.2.

**3. Dimension reduction and algorithms.** In this section we derive algorithms for calculating the local linear estimator with nonadditivity penalty on a grid. In Section 3.2 we derive a formula for computing $\underline{\widehat{\boldsymbol{\beta}}}_R$ which avoids storing and inverting large matrices. An iterative algorithm using these concepts is provided in Section 3.3. Modifications for large $R$ are also discussed.

3.1. *Notation and normal equations.* Define for $k, \kappa = 1, \ldots, d$:

$$S_{0,0}(\underline{\mathbf{x}}) = \frac{1}{n} \sum_{i=1}^n K_h(\underline{\mathbf{X}}_i, \underline{\mathbf{x}}),$$

$$S_{0,k}(\underline{\mathbf{x}}) = S_{k,0}(\underline{\mathbf{x}}) = \frac{1}{n}\sum_{i=1}^{n} K_h(\underline{\mathbf{X}}_i, \underline{\mathbf{x}})\frac{X_{i,k}-x_k}{h_k},$$

$$S_{k,\kappa}(\underline{\mathbf{x}}) = \frac{1}{n}\sum_{i=1}^{n} K_h(\underline{\mathbf{X}}_i, \underline{\mathbf{x}})\frac{X_{i,k}-x_k}{h_k}\frac{X_{i,\kappa}-x_\kappa}{h_\kappa},$$

and for $k = 1, \ldots, d$:

$$L_0(\underline{\mathbf{x}}) = \frac{1}{n}\sum_{i=1}^{n} K_h(\underline{\mathbf{X}}_i, \underline{\mathbf{x}})Y_i,$$

$$L_k(\underline{\mathbf{x}}) = \frac{1}{n}\sum_{i=1}^{n} K_h(\underline{\mathbf{X}}_i, \underline{\mathbf{x}})\frac{X_{i,k}-x_k}{h_k}Y_i.$$

Denote by $\mathbf{S}(\underline{\mathbf{x}})$ the $(d+1) \times (d+1)$ matrix with elements $S_{\ell,l}(\underline{\mathbf{x}})$ and $\underline{\mathbf{L}}(\underline{\mathbf{x}}) = \mathrm{col}_{\ell=0,\ldots,d}(L_\ell(\underline{\mathbf{x}}))$.

Let $\mathbf{S}^{(j)} = \mathbf{S}(\underline{\mathbf{t}}_j)$ and $\underline{\mathbf{L}}^{(j)} = \underline{\mathbf{L}}(\underline{\mathbf{t}}_j)$. The *normal equations* $\widehat{\underline{\boldsymbol{\beta}}}_{ll}^{(j)} = \arg\min_{\underline{\boldsymbol{\beta}}^{(j)}} SSR(\underline{\boldsymbol{\beta}}^{(j)}, \underline{\mathbf{t}}_j)$ for the local linear estimator at $\underline{\mathbf{t}}_j$ are

$$\mathbf{S}^{(j)}\widehat{\underline{\boldsymbol{\beta}}}_{ll}^{(j)} = \underline{\mathbf{L}}^{(j)}.$$

Similarly, for simultaneous local linear estimation on the whole grid we have $\widehat{\underline{\boldsymbol{\beta}}}_{ll} = \mathrm{col}_j(\widehat{\underline{\boldsymbol{\beta}}}_{ll}^{(j)})$, $\mathbf{S} = \mathrm{diag}_j(\mathbf{S}^{(j)})$, $\underline{\mathbf{L}} = \mathrm{col}_j(\underline{\mathbf{L}}^{(j)})$ and

$$\mathbf{S}\widehat{\underline{\boldsymbol{\beta}}}_{ll} = \underline{\mathbf{L}}.$$

The normal equations for the penalized estimator (5) are

$$(6) \qquad (\mathbf{S} + R(\mathbf{I} - \mathbf{P}_{\mathrm{add}}))\widehat{\underline{\boldsymbol{\beta}}}_R = \underline{\mathbf{L}}.$$

3.2. *Dimension reduction.* Simultaneous estimation on a grid requires a large number of parameters. Dimension reduction is necessary for computation.

The normal equations (6) for $\widehat{\underline{\boldsymbol{\beta}}}_R$ are $((\mathbf{S} + R\mathbf{I}) - R\mathbf{P}_{\mathrm{add}})\widehat{\underline{\boldsymbol{\beta}}}_R = \underline{\mathbf{L}}$. Because $\mathbf{S} + R\mathbf{I}$ is a block-diagonal matrix and $R\mathbf{P}_{\mathrm{add}}$ has low rank, solving the normal equations may be simplified using matrix algebra [Rao and Kleffe (1988), page 5, and Appendix A.0.5 here]. We decompose $\mathbf{P}_{\mathrm{add}}$ into a product $\mathbf{Z}^\top\mathbf{Z}$. Using the abbreviation $\mathbf{A}_R = R(\mathbf{S} + R\mathbf{I})^{-1}$, we obtain

$$(7) \qquad \widehat{\underline{\boldsymbol{\beta}}}_R = (\mathbf{I} + \mathbf{A}_R\mathbf{Z}^\top\{\mathbf{I} - \mathbf{Z}\mathbf{A}_R\mathbf{Z}^\top\}^-\mathbf{Z})(\mathbf{S} + R\mathbf{I})^{-1}\underline{\mathbf{L}},$$

where $\{\cdot\}^-$ denotes any generalized inverse.

The matrix $\mathbf{Z}$ has rank $2m^* + 1 - d$, where $m^* = m_1 + \cdots + m_d$. In Appendix A.0.3 an explicit choice for $\mathbf{Z}$ with dimension $2m^* \times m(d+1)$ is given. The multiplication $\mathbf{Z}\underline{\boldsymbol{\beta}}$ consists mainly of $2m^*$ sums of totally $2dm$ terms.

Similarly, calculation of $\mathbf{Z}\mathbf{A}_R\mathbf{Z}^\top$ from $\mathbf{A}_R$ leads to $(2d)^2m$ summations. Calculation of $\mathbf{A}_R$ from $\mathbf{S}$ is of order $d^3m$ operations.

Formula (7) leads to a feasible algorithm because the dimension of the matrices to be inverted is relatively small compared with (6).

*An oblique projection.* Let us define an oblique projection in order to simplify formula (7):

$$\mathbf{P}_{\mathbf{S},R} = \mathbf{Z}^\top\{(\mathbf{I} - \mathbf{Z}\mathbf{Z}^\top) + \mathbf{Z}(\mathbf{I} - \mathbf{A}_R)\mathbf{Z}^\top\}^-\mathbf{Z}(\mathbf{I} - \mathbf{A}_R).$$

In Appendix A.0.5 we show that $\mathbf{P}_{\mathbf{S},R}$ is the orthogonal projection from $\mathbf{F}_{\text{full}}$ to $\mathbf{F}_{\text{add}}$ with respect to the inner product $\langle\underline{\boldsymbol{\beta}},(\mathbf{I} - \mathbf{A}_R)\underline{\boldsymbol{\beta}}\rangle$. In particular, $(\mathbf{I} - \mathbf{A}_R)\mathbf{P}_{\mathbf{S},R}$ is symmetric and $\mathbf{P}_{\mathbf{S},R}^\top(\mathbf{I} - \mathbf{A}_R)(\mathbf{I} - \mathbf{P}_{\mathbf{S},R}) = \mathbf{0}$.

Because $\mathbf{I} - \mathbf{A}_R = (\mathbf{S} + R\mathbf{I})^{-1}\mathbf{S}$ and $\mathbf{S}\underline{\widehat{\boldsymbol{\beta}}}_{ll} = \underline{\mathbf{L}}$, we substitute $(\mathbf{S} + R\mathbf{I})^{-1}\underline{\mathbf{L}}$ in (7) by $(\mathbf{I} - \mathbf{A}_R)\underline{\widehat{\boldsymbol{\beta}}}_{ll}$ and obtain

$$(8) \qquad \underline{\widehat{\boldsymbol{\beta}}}_R = \mathbf{A}_R\mathbf{P}_{\mathbf{S},R}\underline{\widehat{\boldsymbol{\beta}}}_{ll} + (\mathbf{I} - \mathbf{A}_R)\underline{\widehat{\boldsymbol{\beta}}}_{ll}.$$

See Proposition 1 in Section 4.1 for interpretation.

*Modification for large $R$.* For large $R$, $\mathbf{I} - \mathbf{A}_R$ is of order $R^{-1}$ and $\mathbf{P}_{\mathbf{S},R}$ is hence numerically unstable. Because $R(\mathbf{I} - \mathbf{A}_R) = (\mathbf{I} + R^{-1}\mathbf{S})^{-1}\mathbf{S}$ is suitable for large $R$, we modify $\mathbf{P}_{\mathbf{S},R}$ by multiplying both terms $\mathbf{I} - \mathbf{A}_R$ by $R$.

Note that $\mathbf{A}_R = (\mathbf{I} + R^{-1}\mathbf{S})^{-1}$. Formula (8) for large $R$ then becomes

$$\underline{\widehat{\boldsymbol{\beta}}}_R = (R^{-1}\mathbf{I} + \mathbf{A}_R\mathbf{Z}^\top\{(\mathbf{I} - \mathbf{Z}\mathbf{Z}^\top) + \mathbf{Z}(\mathbf{A}_R\mathbf{S})\mathbf{Z}^\top\}^-\mathbf{Z})\mathbf{A}_R\underline{\mathbf{L}}.$$

3.3. *Iterative calculation of the penalized estimator.* We provide in addition an iterative algorithm for the penalized estimator. This avoids inversion of the matrix $\mathbf{I} - \mathbf{Z}\mathbf{A}_R\mathbf{Z}^\top$ and even its calculation.

We use the fact that $\mathbf{P}_{\text{add}}\underline{\widehat{\boldsymbol{\beta}}}_R = \mathbf{P}_{\mathbf{S},R}\underline{\widehat{\boldsymbol{\beta}}}_{ll}$ holds (Proposition 2, Section 4.1) to calculate $\mathbf{P}_{\mathbf{S},R}\underline{\widehat{\boldsymbol{\beta}}}_{ll}$ iteratively via (8),

$$\underline{\widehat{\boldsymbol{\beta}}}_R^{[a+1]} = \mathbf{A}_R\mathbf{P}_{\text{add}}\underline{\widehat{\boldsymbol{\beta}}}_R^{[a]} + (\mathbf{I} - \mathbf{A}_R)\underline{\widehat{\boldsymbol{\beta}}}_{ll}.$$

Only the additive part $\underline{\gamma}^{[a]}$ of $\underline{\widehat{\boldsymbol{\beta}}}_R^{[a]}$ is iterated:

$$(9) \qquad \underline{\gamma}^{[a+1]} = \mathbf{Z}\mathbf{A}_R\mathbf{Z}^\top\underline{\gamma}^{[a]} + \mathbf{Z}(\mathbf{I} - \mathbf{A}_R)\underline{\widehat{\boldsymbol{\beta}}}_{ll},$$

where $\underline{\gamma}^{[a]} = \mathbf{Z}\underline{\widehat{\boldsymbol{\beta}}}_R^{[a]}$. Finally, set

$$\underline{\widehat{\boldsymbol{\beta}}}_R = \mathbf{A}_R\mathbf{Z}^\top\underline{\gamma}^{[\infty]} + (\mathbf{I} - \mathbf{A}_R)\underline{\widehat{\boldsymbol{\beta}}}_{ll}.$$

Uniqueness of (5) implies that $\mathbf{I} - \mathbf{Z}\mathbf{A}_R\mathbf{Z}^\top$ is positive definite [proof in Studer (2002), Appendix B.1]. Accordingly, we have exponential convergence due to fixed point iteration and contraction; see Table 1. In case of nonuniqueness of $\underline{\widehat{\boldsymbol{\beta}}}_R$, the algorithm still converges (see Appendix A.0.6).

The squared difference between the intercepts of $\underline{\widehat{\boldsymbol{\beta}}}_R^{[a]}$ and $\underline{\widehat{\boldsymbol{\beta}}}_R$ in Table 1 diminishes quickly and is negligible for $a \geq 3$ compared with the ISE. The starting value was $\underline{\widehat{\boldsymbol{\beta}}}_R^{[0]} = \underline{\mathbf{0}}$.

*Modification for large $R$.* Algorithm (9) converges because $\mathbf{Z}\mathbf{A}_R\mathbf{Z}^\top$ is a contraction. The eigenvalues of $\mathbf{A}_R$ are, however, increasing with $R$ and $\mathbf{A}_R$ has the identity $\mathbf{I}$ as limit for $R \to \infty$. Therefore convergence is slower for large $R$ and does not work for $R = \infty$.

For large $R$, we choose $\alpha > 0$ such that $\alpha\mathbf{S} < \mathbf{I}$ and use

$$(10) \qquad \mathbf{Z}\underline{\widehat{\boldsymbol{\beta}}}_R = \mathbf{Z}(\mathbf{I} - \alpha R(\mathbf{I} - \mathbf{A}_R))\mathbf{Z}^\top\mathbf{Z}\underline{\widehat{\boldsymbol{\beta}}}_R + \alpha\mathbf{Z}R(\mathbf{I} - \mathbf{A}_R)\underline{\widehat{\boldsymbol{\beta}}}_{ll}$$

for iterations instead of (9).

*Generalizations.* The derivations for (7) assume only that $\mathbf{Z}^\top\mathbf{Z}$ is a projection. Hence, if $\mathbf{F}_{\mathrm{add}}$ is replaced by another subspace $\mathbf{F}_{\mathrm{sub}}$, say, and $\mathbf{Z}$ is modified such that $\mathbf{Z}^\top\mathbf{Z}$ is the orthogonal projection from $\mathbf{F}_{\mathrm{full}}$ to $\mathbf{F}_{\mathrm{sub}}$, then the above algorithms remain valid. Generalization from local linear to local polynomial estimation is achieved by corresponding modification of $\mathbf{S}$ and $\underline{\mathbf{L}}$.

Implementation is simplified by the fact that $\mathbf{Z}$ need not have full rank. For the iterative algorithm (9), moreover, there is no need to calculate the matrix $\mathbf{Z}$ explicitly. For example, if $\mathbf{F}_{\mathrm{sub}}$ corresponds to using bivariate interaction terms in the additive model or postulating the same regression function for subgroups, multiplication by $\mathbf{Z}$, $\mathbf{Z}^\top$ and $\mathbf{Z}\mathbf{A}_R\mathbf{Z}^\top$ can be implemented efficiently.

**4. Properties of the estimator.** In this section, we evaluate the effect of the nonadditivity penalty on the estimator. Both on a grid (Section 4.1) and on an interval (Section 4.2), the penalized estimator turns out to be a pointwise compromise between the local linear and some ($R$-dependent) additive estimator. The compromise depends on how well the local linear estimator is determined locally. This is an attractive property as it leads to automatic regularization in sparse regions (provided that the additive estimator is well determined). The additive part of $\widehat{\underline{r}}_R$ is studied using two different norms (Propositions 3 and 4).

TABLE 1
*Convergence of iterations for the estimator $\widehat{r}_R$ in Figure 2(d)*

| $a$ | 0 | 1 | 2 | 3 | 4 | 5 | ISE |
|---|---|---|---|---|---|---|---|
| $\|[\underline{\widehat{\boldsymbol{\beta}}}_R^{[a]} - \underline{\widehat{\boldsymbol{\beta}}}_R]_{\mathrm{intercept}}\|^2$ | 463 | 15.5 | 0.8 | 0.1 | 0.02 | 0.008 | 6.0 |

Later on, we focus on the smoothness of the model choice via the penalty parameter $R$ for fixed $n$. Continuity in $R$ for $R \in (0, \infty)$ is obvious and the cases $R = 0$ and $\infty$ are investigated in Section 4.4. We investigate the rate of convergence of $\widehat{\underline{r}}_R$ to $\widehat{\underline{r}}_{\text{add}}$, depending on whether or not $r^{\text{true}}$ is additive. In both cases we find a rate for $R$ such that $\| \widehat{\underline{r}}_R - \widehat{\underline{r}}_{\text{add}} \|_2^2$ is of smaller order than $n^{-4/5}$. In Section 4.5 we consider the data-adaptive choice of $R$ and $h$. In Section 4.6 we see that in the case of fixed uniform design ($d \leq 4$) $\widehat{\underline{r}}_R$ with data-adaptive $R$ is equivalent to $\widehat{\underline{r}}_{\text{add}}$ for additive functions.

4.1. *Properties of the estimator on a grid.* We investigate the effect of the penalty on estimation at one output point: the penalized estimator is a kind of convex combination between a local linear and an additive estimator. Furthermore, the local linear estimator may be decomposed into a sum of additive and residual components. The penalized estimator is the sum of the additive part and shrunken residuals, which are orthogonal to the additive part.

In Section 3.2 an oblique projection $\mathbf{P}_{\mathbf{S},R}$ was introduced, leading to

$$\widehat{\underline{\beta}}_R = \mathbf{A}_R \mathbf{P}_{\mathbf{S},R} \widehat{\underline{\beta}}_{ll} + (\mathbf{I} - \mathbf{A}_R) \widehat{\underline{\beta}}_{ll}$$

in (8). Recall that $\mathbf{A}_R = \text{diag}_j(R(\mathbf{S}^{(j)} + R\mathbf{I})^{-1})$ is *block-diagonal* with eigenvalues between zero and 1. Let us see what (8) implies for one output point. Denote by $(\mathbf{P}_{\mathbf{S},R} \widehat{\underline{\beta}}_{ll})^{(j)}$ the components of $\mathbf{P}_{\mathbf{S},R} \widehat{\underline{\beta}}_{ll}$ corresponding to output point $\underline{t}_j$, formally $\mathbf{P}_{\mathbf{S},R} \widehat{\underline{\beta}}_{ll} = \text{col}_j((\mathbf{P}_{\mathbf{S},R} \widehat{\underline{\beta}}_{ll})^{(j)})$.

PROPOSITION 1. *The penalized estimator $\widehat{\underline{\beta}}_R^{(j)}$ is a pointwise compromise between some ($R$-dependent) additive fit $(\mathbf{P}_{\mathbf{S},R} \widehat{\underline{\beta}}_{ll})^{(j)}$ and the local linear fit $\widehat{\underline{\beta}}_{ll}^{(j)}$:*

$$\widehat{\underline{\beta}}_R^{(j)} = (\mathbf{S}^{(j)} + R\mathbf{I})^{-1} \{ R(\mathbf{P}_{\mathbf{S},R} \widehat{\underline{\beta}}_{ll})^{(j)} + \mathbf{S}^{(j)} \widehat{\underline{\beta}}_{ll}^{(j)} \}.$$

In sparse regions the local linear estimator is unstable [Seifert and Gasser (1996)], because $\mathbf{S}^{(j)}$ may be nearly singular. The above formula indicates that penalizing solves this problem as a byproduct, because the additive part of $\widehat{\underline{\beta}}_R^{(j)}$ is stable under weaker conditions. This regularization property is illustrated in Figure 2(b) versus 2(d). When all eigenvalues of $\mathbf{S}^{(j)}$ are large, the effect of a small penalty $R$ vanishes.

We derive now a decomposition of $\widehat{\underline{\beta}}_R$ into an additive component and an orthogonal remainder term. Formula (8) is equivalent to

(11) $$\widehat{\underline{\beta}}_R = \mathbf{P}_{\mathbf{S},R} \widehat{\underline{\beta}}_{ll} + (\mathbf{I} - \mathbf{A}_R)(\mathbf{I} - \mathbf{P}_{\mathbf{S},R}) \widehat{\underline{\beta}}_{ll}.$$

Only the nonadditive part involves shrinking.

PROPOSITION 2. *The following relations hold:*

$$\mathbf{P}_{\mathrm{add}}\widehat{\underline{\boldsymbol{\beta}}}_R = \mathbf{P}_{\mathbf{S},R}\widehat{\underline{\boldsymbol{\beta}}}_{ll} \quad and \quad (\mathbf{I}-\mathbf{P}_{\mathrm{add}})\widehat{\underline{\boldsymbol{\beta}}}_R = (\mathbf{I}-\mathbf{A}_R)(\mathbf{I}-\mathbf{P}_{\mathbf{S},R})\widehat{\underline{\boldsymbol{\beta}}}_{ll}.$$

The proof is in Appendix A.0.7.

4.2. *Properties of the estimator on an interval.* Now we will show that Propositions 1 and 2 hold not only on a grid but also on an interval. Proposition 4 states that the additive part of $\widehat{\underline{r}}_R$ with respect to $\mathcal{P}_*$ is $\widehat{\underline{r}}_{\mathrm{add}}$, independent of $R$.

Define the symmetric, continuous operator $\mathcal{S}_* : \mathcal{F}_{\mathrm{full}} \to \mathcal{F}_{\mathrm{full}}$, $\underline{r} \mapsto \breve{\underline{r}}$ by

$$\begin{pmatrix} \breve{r}^0(\underline{\mathbf{x}}) \\ \vdots \\ \breve{r}^d(\underline{\mathbf{x}}) \end{pmatrix} = \begin{pmatrix} S_{0,0}(\underline{\mathbf{x}})r^0(\underline{\mathbf{x}}) & + & \cdots & + & S_{0,d}(\underline{\mathbf{x}})r^d(\underline{\mathbf{x}}) \\ & \vdots & & & \vdots \\ S_{d,0}(\underline{\mathbf{x}})r^0(\underline{\mathbf{x}}) & + & \cdots & + & S_{d,d}(\underline{\mathbf{x}})r^d(\underline{\mathbf{x}}) \end{pmatrix}$$

(see Section 3.1). We have by construction that $\|\underline{r}\|_*^2 = \langle \underline{r}, \mathcal{S}_*\underline{r} \rangle_2$. Let $\underline{r}_L \in \mathcal{F}_{\mathrm{full}}$ with $r_L^\ell(\underline{\mathbf{x}}) = L_\ell(\underline{\mathbf{x}})$. The *normal equations* for $\widehat{\underline{r}}_{ll}$, the minimizer of (3) in Section 2.1, are

$$\mathcal{S}_*\widehat{\underline{r}}_{ll} = \underline{r}_L.$$

The *normal equations* for $\widehat{\underline{r}}_R$, the minimizer of (4), are

$$(12) \qquad (\mathcal{S}_* + R(\mathcal{I}-\mathcal{P}_{\mathrm{add}}))\widehat{\underline{r}}_R = \underline{r}_L \equiv \mathcal{S}_*\widehat{\underline{r}}_{ll}.$$

Let $\mathcal{P}_{*,R}$ denote the orthogonal projection from $\mathcal{F}_{\mathrm{full}}$ to $\mathcal{F}_{\mathrm{add}}$ with respect to the norm $\|\cdot\|_R$. According to MLN, $\mathcal{P}_{*,R}$ is continuous with probability tending to 1 for $n \to \infty$, under regularity conditions for the design density and kernel [Conditions MLN:B1 and MLN:B2$'$ in Appendix A.0.2]. In particular, the bandwidth $h$ is of order $n^{-1/5}$ or larger [Condition C1+]. An explicit formula for $\mathcal{P}_{*,R}$ is given in (26), Appendix A.0.8.

Then $\widehat{\underline{r}}_R$ may be decomposed similarly to (8) and Proposition 1:

$$(13) \qquad \widehat{\underline{r}}_R = \{(\mathcal{S}_* + R\mathcal{I})^{-1}R\mathcal{P}_{*,R} + (\mathcal{S}_* + R\mathcal{I})^{-1}\mathcal{S}_*\}\widehat{\underline{r}}_{ll}.$$

Note that $(\mathcal{S}_* + R\mathcal{I})^{-1}\mathcal{S}_*$ is a *pointwise* (in $\underline{\mathbf{x}}$) matrix multiplication. Furthermore, $(\mathcal{S}_* + R\mathcal{I})^{-1}\mathcal{S}_*$ and $(\mathcal{S}_* + R\mathcal{I})^{-1}R$ sum to $\mathcal{I}$ and have eigenvalues between zero and 1. Hence, (13) indicates that $\widehat{\underline{r}}_R(\underline{\mathbf{x}})$ is *some kind of convex combination* of $\widehat{\underline{r}}_{ll}$ and $\mathcal{P}_{*,R}\widehat{\underline{r}}_{ll}$.

Similarly to Proposition 2, the above formula may be rewritten as

$$(14) \qquad \widehat{\underline{r}}_R = \{\mathcal{P}_{*,R} + (\mathcal{S}_* + R\mathcal{I})^{-1}\mathcal{S}_*(\mathcal{I}-\mathcal{P}_{*,R})\}\widehat{\underline{r}}_{ll}.$$

PROPOSITION 3.  *The following relations hold:*

$$\mathcal{P}_{\mathrm{add}} \widehat{\underline{r}}_R = \mathcal{P}_{*,R} \widehat{\underline{r}}_{ll} \quad and \quad (\mathfrak{I} - \mathcal{P}_{\mathrm{add}}) \widehat{\underline{r}}_R = (\mathcal{S}_* + R\mathfrak{I})^{-1} \mathcal{S}_* (\mathfrak{I} - \mathcal{P}_{*,R}) \widehat{\underline{r}}_{ll}.$$

In Section 4.4 we will see that $\mathcal{P}_{\mathrm{add}} \widehat{\underline{r}}_R - \widehat{\underline{r}}_{\mathrm{add}}$ is $O(R^{-1})$, for fixed $n$. The $R$-dependence of the additive part $\mathcal{P}_{\mathrm{add}} \widehat{\underline{r}}_R$ can be avoided when using the oblique projection $\mathcal{P}_*$ instead:

PROPOSITION 4.  $\mathcal{P}_* \widehat{\underline{r}}_R = \widehat{\underline{r}}_{\mathrm{add}}$ *holds.*

This means that the MLN estimator $\widehat{\underline{r}}_{\mathrm{add}}$ is the additive part of $\widehat{\underline{r}}_R$ with respect to the norm $\| \cdot \|_*$. Both proofs can be found in Appendix A.0.8.

4.3. *Bounding $\mathcal{S}_*$ and $\mathcal{P}_{*,R}$.*  If $S_{0,0}(\underline{\mathbf{x}})$ is a uniformly consistent density estimator, the operators $\mathcal{S}_*$ and $\mathcal{P}_{*,R}$ are shown to be bounded. This property will be used in Section 4.4.

Let $\|\mathcal{S}_*\|_{2,\mathrm{sup}}$ denote the supremum norm of $\mathcal{S}_*$ based on the Euclidean norm on $\mathcal{F}_{\mathrm{full}}$. Here, we want to find upper bounds for $\|\mathcal{S}_*\|_{2,\mathrm{sup}}$, that is, a uniform bound for the maximum eigenvalue of $\mathbf{S}(\underline{\mathbf{x}})$. Because the kernel is bounded with compact support by Condition MLN:B1, the maximal eigenvalue of $\mathbf{S}(\underline{\mathbf{x}})$ is of order $S_{0,0}(\underline{\mathbf{x}})$. Note that $S_{0,0}(\underline{\mathbf{x}})$ is a kernel density estimator of $f$. We are thus interested in uniform boundedness from above of $\widehat{f}(\underline{\mathbf{x}}) = S_{0,0}(\underline{\mathbf{x}})$.

Silverman (1978) derived uniform consistency of kernel density estimators for $d = 1$. We will use a result of Gao (2003), which asserts uniform consistency for density estimators for *continuous densities* on $\mathbb{R}^d$ and bandwidths $h$ satisfying

(15)           $h \to 0$  and  $\dfrac{nh^d}{\log(h^{-1})} \to \infty$     as $n \to \infty$.

By Condition C1+ these conditions are satisfied for $d \le 4$. For $d \ge 5$, the condition (15) is not satisfied for the optimal bandwidth $\propto n^{-1/5}$ of the additive model. We thus lose flexibility in the model choice when (15) is assumed.

PROPOSITION 5.  *Under Conditions* MLN:B1, MLN:B2$'$ *and* (15), $\|\mathcal{S}_*\|_{2,\mathrm{sup}}$ *is uniformly bounded with probability tending to 1 for $n \to \infty$.*

Note that for fixed $n$, this norm is always bounded because of Condition MLN:B1. The $R$-dependent projection $\mathcal{P}_{*,R}$ may be bounded uniformly in $R$ using Proposition 5 and Lemma 2 in Section 4.4:

LEMMA 1.  *Under Conditions* MLN:B1, MLN:B2$'$ *and* C1+$\cap$ (15), $\|\mathcal{P}_{*,R}\|_{2,\mathrm{sup}} = O_P(1)$, *uniformly in $R$.*

The proofs are in Appendix A.0.9.

4.4. *The additive and full models as special cases.* Here we justify the interpretation of $\widehat{\underline{r}}_{ll}$ and $\widehat{\underline{r}}_{\mathrm{add}}$ as special cases of $\widehat{\underline{r}}_R$ for $R = 0$ and $R = \infty$, respectively. This is appreciated because $\widehat{\underline{r}}_{ll}$ and $\widehat{\underline{r}}_{\mathrm{add}}$ are known to be asymptotically optimal for the respective situations. The rate of convergence of $\widehat{\underline{r}}_R$ to $\widehat{\underline{r}}_{\mathrm{add}}$ for $R \to \infty$ depends on the supremum norm of $\mathcal{S}_*$ and whether or not the regression function is additive.

We will start with the convergence of $\widehat{\underline{r}}_R$ to $\widehat{\underline{r}}_{ll}$ for $R \downarrow 0$. Consider the case where $\mathbf{S}(\underline{\mathbf{x}})^{-1}$ is uniformly continuous in $\underline{\mathbf{x}}$. This is a sufficient condition for *bounded variance* of $\widehat{\underline{r}}_{ll}$ and represents therefore the *well-behaved cases*. In this case, $\mathcal{S}_*$ has a continuous inverse and the limit of $\widehat{\underline{r}}_R$ for $R \downarrow 0$ is $\widehat{\underline{r}}_{ll}$. Let us mention that uniform continuity is a stronger assumption than uniqueness of $\widehat{\underline{r}}_{ll}$. Uniform continuity means that for any $\underline{r} \in \mathcal{F}_{\mathrm{full}}$ with $\|\underline{r}\|_2 = 1$ the norm $\|\underline{r}\|_*$ is *bounded away from zero*, whereas uniqueness needs only a *nonzero* norm. If $\widehat{\underline{r}}_{ll}$ is not well determined, a small positive penalty provides the desired regularization.

[Mammen, Linton and Nielsen](1999) showed that the *additive* estimator

$$\widehat{\underline{r}}_{\mathrm{add}} = \arg\min_{\underline{r} \in \mathcal{F}_{\mathrm{add}}} \|\underline{r}_Y - \underline{r}\|_*^2$$

is *asymptotically oracle optimal* under Conditions MLN:B1–B4′ and C1 for additive regression functions as in (1). Let us therefore examine the convergence of $\widehat{\underline{r}}_R$ to $\widehat{\underline{r}}_{\mathrm{add}}$ for $R \to \infty$. Decompose via

(16) $$\|\widehat{\underline{r}}_R - \widehat{\underline{r}}_{\mathrm{add}}\|_2^2 = \|(\mathcal{I} - \mathcal{P}_{\mathrm{add}})\widehat{\underline{r}}_R\|_2^2 + \|\mathcal{P}_{\mathrm{add}}\widehat{\underline{r}}_R - \widehat{\underline{r}}_{\mathrm{add}}\|_2^2.$$

For bounds of the first term of the sum, see Lemma 4 below.

Recall that $\mathcal{P}_{\mathrm{add}}\widehat{\underline{r}}_R = \mathcal{P}_{*,R}\widehat{\underline{r}}_{ll}$ holds by Proposition 3. Similar to the modifications for large $R$ in Section 3.2, we introduce an alternative formula for $\mathcal{P}_{*,R}\widehat{\underline{r}}_{ll}$ [defined in (26), Appendix A.0.8],

(17) $$\mathcal{P}_{\mathrm{add}}\widehat{\underline{r}}_R = (\mathcal{P}_{\mathrm{add}}(\mathcal{I} + R^{-1}\mathcal{S}_*)^{-1}\mathcal{S}_*\mathcal{P}_{\mathrm{add}})_{|\mathcal{F}_{\mathrm{add}}}^{-1}\mathcal{P}_{\mathrm{add}}(\mathcal{I} + R^{-1}\mathcal{S}_*)^{-1}\underline{r}_L,$$

where $(\cdots)_{|\mathcal{F}_{\mathrm{add}}}$ indicates that the expression is an operator on $\mathcal{F}_{\mathrm{add}}$. Define $\mathcal{S}_{\mathrm{add}} = (\mathcal{P}_{\mathrm{add}}\mathcal{S}_*\mathcal{P}_{\mathrm{add}})_{|\mathcal{F}_{\mathrm{add}}}$. Solving the normal equations for $\widehat{\underline{r}}_{\mathrm{add}}$ leads to

$$\widehat{\underline{r}}_{\mathrm{add}} = \mathcal{S}_{\mathrm{add}}^{-1}\mathcal{P}_{\mathrm{add}}\underline{r}_L.$$

The right-hand side is equal to the right-hand side of $\mathcal{P}_{\mathrm{add}}\widehat{\underline{r}}_R$ for $R^{-1} = 0$. If $R^{-1}\|\mathcal{S}_*\|_{2,\mathrm{sup}}$ tends to zero, we may use a Taylor approximation and obtain $\|\mathcal{P}_{\mathrm{add}}\widehat{\underline{r}}_R - \widehat{\underline{r}}_{\mathrm{add}}\|_2 = O_P(R^{-1}\|\mathcal{S}_*\|_{2,\mathrm{sup}})$; see Lemmas 2 and 3.

LEMMA 2. *Under Conditions* MLN:B1, MLN:B2′ *and* C1+, $\mathcal{S}_{\mathrm{add}}$ *converges for* $n \to \infty$ *to an operator with continuous inverse.*

*Hence,* $\mathcal{S}_{\mathrm{add}}^{-1}$ *is continuous with probability tending to 1 for* $n \to \infty$ *and* $\|\mathcal{S}_{\mathrm{add}}^{-1}\|_{2,\mathrm{sup}}$ *is uniformly bounded in* $n$, $\forall n \geq \tilde{n}$, *with probability tending to 1 for* $\tilde{n} \to \infty$.

The proof is given in Appendix A.0.10. Uniqueness of $\widehat{\underline{r}}_{\text{add}}$ is equivalent to $\|\underline{r}_{\text{add}}\|_* > 0$ for all $\underline{r}_{\text{add}} \in \mathcal{F}_{\text{add}} - \{\underline{0}\}$. Lemma 2 states that $\|\underline{r}_{\text{add}}\|_*/\|\underline{r}_{\text{add}}\|_2$ is bounded away from zero with probability tending to 1.

LEMMA 3. *Under Conditions* MLN:B1, MLN:B2′, MLN:B3′ *and* C1+ *we obtain for fixed $n$ (i.e., conditional on the data)*

$$\sup_{\underline{\mathbf{x}}} |\mathcal{P}_{\text{add}}(\mathcal{I} - (\mathcal{I} + R^{-1}\mathcal{S}_*)^{-1})\underline{r}_L(\underline{\mathbf{x}})| = O(R^{-1})$$

*and* $\|\mathcal{P}_{\text{add}}\underline{r}_L\|_2$ *is finite.*

*For increasing $n$, we obtain*

$$\sup_{\underline{\mathbf{x}}} |\mathcal{P}_{\text{add}}(\mathcal{I} - (\mathcal{I} + R^{-1}\mathcal{S}_*)^{-1})\underline{r}_L(\underline{\mathbf{x}})| = O_P(R^{-1}\|\mathcal{S}_*\|_{2,\text{sup}}).$$

*Furthermore,*

$$\|\mathcal{P}_{\text{add}}\underline{r}_L\|_2 = O_P(1).$$

A proof is given in Appendix A.0.10.

LEMMA 4. *The following bound holds:*

$$\|(\mathcal{I} - \mathcal{P}_{\text{add}})\widehat{\underline{r}}_R\|_2^2 \leq 2R^{-1}\|\mathcal{S}_*\|_{2,\text{sup}}\|\underline{r}_Y - \widehat{\underline{r}}_{\text{add}}\|_*\|\widehat{\underline{r}}_R - \widehat{\underline{r}}_{\text{add}}\|_2.$$

*Under Condition* MLN:B3′ *we have*

$$\|\underline{r}_Y - \widehat{\underline{r}}_{\text{add}}\|_*^2 \leq \frac{1}{n}\sum_{i=1}^n Y_i^2 = \begin{cases} \text{finite,} & \text{for fixed } n \text{ (and } \underline{\mathbf{Y}}), \\ O_P(1), & \text{for increasing } n. \end{cases}$$

See Appendix A.0.10 for a proof. Using (16) we obtain:

THEOREM 1. *Assume Conditions* MLN:B1, MLN:B2′, MLN:B3′ *and* C1+. *For fixed $n$,*

$$\|\widehat{\underline{r}}_R - \widehat{\underline{r}}_{\text{add}}\|_2 = O(R^{-1})$$

*holds with probability tending to 1 for increasing $n$. Formally, this means* $P[\limsup_{R\to\infty} R\|\widehat{\underline{r}}_R - \widehat{\underline{r}}_{\text{add}}\|_2 < \infty] \overset{n\to\infty}{\to} 1.$

*For $n \to \infty$*

$$\|\widehat{\underline{r}}_R - \widehat{\underline{r}}_{\text{add}}\|_2 = O_P(R^{-1}\|\mathcal{S}_*\|_{2,\text{sup}}).$$

Note that this holds also for *nonadditive* regression functions. For additive regression functions we obtain a better bound; see Theorem 2.

Applying Proposition 5 to Theorem 1, for $h \propto n^{-1/5}$ (Condition C1), $d \leq 4$ and $R^{-1} = o(n^{-2/5})$, we get $\|\widehat{\underline{r}}_R - \widehat{\underline{r}}_{\text{add}}\|_2^2 = o_P(n^{-4/5})$. For $h \propto n^{-1/5}$ and $d \geq 5$, $\|\mathcal{S}_*\|_{2,\text{sup}}$ is not bounded by a constant and $R$ needs to converge faster to $\infty$ to achieve equivalence. Alternatively, one might use a larger bandwidth. (Without proof.)

THEOREM 2. *Assume Conditions* MLN:B1–B4′, $d \leq 4$, $h \propto [n^{-1/5}, n^{-1/(4+d)}]$ *and an* additive *regression function* (1). *Then*

$$\|\underline{\hat{r}}_R - \underline{\hat{r}}_{\text{add}}\|_2 = O_P\left(\frac{1}{R\sqrt{nh^d}}\right).$$

For $h \propto n^{-1/5}$ and additive regression functions, $R^{-1} = o(n^{-(d-1)/10})$ is sufficient to obtain equivalence of $\underline{\hat{r}}_R$ and $\underline{\hat{r}}_{\text{add}}$. The proof is in Appendix A.0.10.

4.5. *Data-adaptive parameter selection.* We consider the simultaneous choice of the tuning parameters $R$ and $h$. In the case of an additive regression function $r^{\text{true}}$, the first-order bias of $\underline{\hat{r}}_R$ is independent of $R$ and parameter selection is asymptotically equivalent to the classical variance/bias compromise. Hence, $\hat{h} \propto n^{-1/5}$ and $\hat{R} \to \infty$. The rate of $\hat{R}$ is investigated in Section 4.6.

Asymptotically, we have only to consider the cases $r^{\text{true}} =$ additive or full model, and the question is then whether $\hat{R}$ is able to identify these cases.

We consider parameter selection criteria that depend on fitted values at design points. The vector of fitted values $\underline{\hat{Y}} = \text{col}_{i=1,\ldots,n}(\hat{r}^0_R(\underline{X}_i)) = \mathbf{M}_R \underline{Y}$ *depends linearly* on $\underline{Y}$, where $\mathbf{M}_R$ is called "hat matrix."

In practice, the estimator is computed on a grid and we need some *interpolation* to obtain estimates at the design points. Define

$$\underline{\hat{\boldsymbol{\beta}}}_R(\underline{X}_i) = (\mathbf{S}(\underline{X}_i) + R\mathbf{I})^{-1}\{\mathbf{S}(\underline{X}_i)\underline{\hat{\boldsymbol{\beta}}}_{ll}(\underline{X}_i) + R(\mathbf{P}_{\mathbf{S},R}\underline{\hat{\boldsymbol{\beta}}}_{ll})(\underline{X}_i)\},$$

where $(\mathbf{P}_{\mathbf{S},R}\underline{\hat{\boldsymbol{\beta}}}_{ll})(\underline{X}_i)$ denotes the interpolated value at $\underline{X}_i$ of the additive part $\mathbf{P}_{\mathbf{S},R}\underline{\hat{\boldsymbol{\beta}}}_{ll}$. Therefore, $\underline{\hat{Y}} = \text{col}_i([\underline{\hat{\boldsymbol{\beta}}}_R(\underline{X}_i)]_0)$ is a linear combination of $\underline{Y}$ and construction of $\mathbf{M}_R$ is obvious.

We consider the following criteria:

$$\text{AIC}(R, h) = \log(\hat{\sigma}^2) + 2\text{tr}(\mathbf{M}_R)/n,$$

$$\text{GCV}(R, h) = \frac{\hat{\sigma}^2}{(1 - \text{tr}(\mathbf{M}_R)/n)^2},$$

$$\text{AIC}_{\text{C}}(R, h) = \log(\hat{\sigma}^2) + \frac{1 + \text{tr}(\mathbf{M}_R)/n}{1 - (\text{tr}(\mathbf{M}_R) + 2)/n},$$

where $\hat{\sigma}^2 = \frac{1}{n}\|\underline{Y} - \mathbf{M}_R\underline{Y}\|^2$ and $\text{tr}(\mathbf{M}_R)$ denotes the *trace* of $\mathbf{M}_R$, which is interpreted as *degrees of freedom*. AIC and GCV are classical model selection criteria [see, e.g., Hastie and Tibshirani (1990)], and AIC$_{\text{C}}$ was introduced by Hurvich, Simonoff and Tsai (1998). These criteria are justified for $\underline{\hat{r}}_{ll}$ only when (15) is satisfied. As we want to analyze the ability of $\hat{R}$ to identify the additive model with its optimal bandwidth (Condition C1), we will assume $d \leq 4$ throughout this section. Moreover, we will assume that $\mathbb{E}[\varepsilon^4] < \infty$.

Let us compare the criteria $(\tau = \frac{1}{n}\mathrm{tr}(\mathbf{M}_R))$

$$\mathrm{AIC} = \log(\widehat{\sigma}^2) + 2\tau,$$

$$\log(\mathrm{GCV}) = \log(\widehat{\sigma}^2) + 2\left(\tau + \frac{\tau^2}{2} + \frac{\tau^3}{3} + \cdots\right),$$

$$\mathrm{AIC_C} - 1 = \log(\widehat{\sigma}^2) + \frac{2}{n-2} + 2\sum_{k \geq 1} \frac{(n-1)n^k \tau^k}{(n-2)^{k+1}}$$

$$\approx \log(\widehat{\sigma}^2) + 2(\tau + \tau^2 + \tau^3 + \cdots).$$

All are of the form $\log(\widehat{\sigma}^2)$ plus some penalty against undersmoothing; see Härdle, Hall and Marron (1988). As the penalty increases from AIC to $\log(\mathrm{GCV})$ and further to $\mathrm{AIC_C}$, minimizing these criteria leads to increasingly more smoothing (decreasing $\tau_{\min}$): According to Hurvich, Simonoff and Tsai, $\mathrm{AIC_C}$ avoids the large variability and the tendency to undersmooth of GCV and classical AIC observed when estimating *bandwidths* for $d = 1$. Note that AIC has its global minimum at interpolation $(h = 0, R = 0)$, leading to $\widehat{\sigma}^2 = 0$ and $\tau = 1$. Undersmoothing, however, contradicts the aim of this paper and is avoided by assuming $h \propto [n^{-1/5}, n^{-1/(4+d)}]$ and $R \geq R_{\min}(h)$. The lower bound $R_{\min}(h)$ is chosen to bound the variance of $\widehat{\underline{r}}_R$ by the optimum rate $n^{-4/(4+d)}$. This condition does not rule out the asymptotically optimal additive estimator $(R = \infty, h \propto n^{-1/5})$. Then $\frac{1}{n}\mathrm{tr}(\mathbf{M}_R) \to 0$ as $n \to \infty$ and $\frac{\widehat{\sigma}^2}{\sigma^2} - 1$ is $O_P(\frac{1}{\sqrt{n}})$.

Hence, we may use the approximation $\log(\widehat{\sigma}^2) = \log(\sigma^2) + \frac{\widehat{\sigma}^2}{\sigma^2} - 1 + O_P(\frac{1}{n})$. Define the Taylor approximation of $\mathrm{AIC} - \log(\sigma^2)$,

$$\mathrm{AIC_T} = \frac{\widehat{\sigma}^2}{\sigma^2} - 1 + \frac{2}{n}\mathrm{tr}(\mathbf{M}_R).$$

The expected value of $\widehat{\underline{r}}_{ll}$ is for additive regression functions

$$\mathbb{E}(\widehat{\underline{r}}_{ll}) = \left(r_{\mathrm{add}}^{\mathrm{true},0}(\underline{\mathbf{x}}) + \frac{\mu_2(K)}{2}\sum_{k=1}^{d} h_k^2 \frac{\partial^2}{\partial x_k^2} r_{\mathrm{add},k}^{\mathrm{true}}(x_k),\right.$$

(18)

$$\left. h_1 \frac{\partial}{\partial x_1} r_{\mathrm{add},1}^{\mathrm{true}}(x_1), \ldots, h_d \frac{\partial}{\partial x_d} r_{\mathrm{add},d}^{\mathrm{true}}(x_d)\right) + o_p(h^2).$$

The leading terms are in $\mathcal{F}_{\mathrm{add}}$ and hence unchanged under multiplication by $\mathcal{P}_{*,R}$, and the $o_P(h^2)$ terms remain small enough because of Lemma 1. Accordingly, $(\mathbf{I} - \mathbf{M}_R)\underline{r}_{\mathrm{add}}^{\mathrm{true}} = O(h^2 \underline{\mathbf{1}})$, and the first-order terms of the bias of $\widehat{\underline{\mathbf{Y}}}$ are independent of $R$.

Next, we need to ensure that $\widehat{\underline{r}}_R$ is not degenerate. Using (14) in Section 4.2, we choose some small constant $R_{\min} > 0$, assume that $R \geq R_{\min}$, and consequently $(\mathcal{S}_* + R\mathcal{I})^{-1}\mathcal{S}_* \widehat{\underline{r}}_{ll}$ is stable (ridge regression). With probability tending to 1, $\mathcal{S}_{\mathrm{add}}^{-1}$ is continuous (Lemma 2), that is, $\widehat{\underline{r}}_{\mathrm{add}}$ is stable. If both $\mathcal{S}_*$ and $\mathcal{S}_{\mathrm{add}}^{-1}$ are bounded, we do not have to worry about stability of $\mathcal{P}_{*,R}\widehat{\underline{r}}_{ll}$ (see also the proof of Lemma 1). Therefore, $\|\mathbf{M}_R^\top \mathbf{M}_R\|_{\sup} = O_P(1)$. Obviously, when $\mathcal{S}_*^{-1}$ is continuous, we need not assume that $R \geq R_{\min}$.

Therefore, $\frac{1}{n\sigma^2}\langle(\mathbf{I} - \mathbf{M}_R)\underline{\varepsilon}, (\mathbf{I} - \mathbf{M}_R)\underline{r}_{\mathrm{add}}^{\mathrm{true}}\rangle = O_P(\frac{h^2}{\sqrt{n}})$, which is $o_P(h^4)$ for $h$ as in Condition C1+. Note that $\sigma^2 \operatorname{tr}(\mathbf{M}_R) = \mathbb{E}[\langle\underline{\varepsilon}, \mathbf{M}_R\underline{\varepsilon}\rangle]$. Hence

$$\mathrm{AIC_T} - \left(\frac{1}{n\sigma^2}\|\underline{\varepsilon}\|^2 - 1\right)$$
$$= \frac{1}{n\sigma^2}\mathbb{E}[\|\mathbf{M}_R\underline{\varepsilon}\|^2] + \frac{1}{n\sigma^2}\|(\mathbf{I} - \mathbf{M}_R)\underline{r}_{\mathrm{add}}^{\mathrm{true}}\|^2$$
$$+ \frac{1}{n\sigma^2}(1 - \mathbb{E})[\|\mathbf{M}_R\underline{\varepsilon}\|^2 + 2\langle\underline{\varepsilon}, \mathbf{M}_R\underline{\varepsilon}\rangle] + O_P\left(\frac{h^2}{\sqrt{n}}\right),$$

where $(1 - \mathbb{E})[\langle\underline{\varepsilon}, \mathbf{M}\underline{\varepsilon}\rangle] = \langle\underline{\varepsilon}, \mathbf{M}\underline{\varepsilon}\rangle - \mathbb{E}[\langle\underline{\varepsilon}, \mathbf{M}\underline{\varepsilon}\rangle]$.

LEMMA 5. *For $\mathbb{E}[\varepsilon^4] < \infty$, $\operatorname{var}(\langle\underline{\varepsilon}, \mathbf{M}_R\underline{\varepsilon}\rangle) = O(\mathbb{E}[\|\mathbf{M}_R\underline{\varepsilon}\|^2])$. Moreover, if $\|\mathbf{M}_R^\top \mathbf{M}_R\|_{\sup} = O_P(1)$, then $\operatorname{var}(\|\mathbf{M}_R\underline{\varepsilon}\|^2) = O_P(\mathbb{E}[\|\mathbf{M}_R\underline{\varepsilon}\|^2])$.*

Because $\mathbb{E}(\|\mathbf{M}_R\underline{\varepsilon}\|^2)$ is of order $h^{-d}$ and $h^{-1}$ for $R = 0$ and $R = \infty$, respectively, the standard deviation of $(1 - \mathbb{E})[\cdots]$ is of smaller order. The proof is given in Appendix A.0.11.

The leading terms of $\mathrm{AIC_T} - (\frac{1}{n\sigma^2}\|\underline{\varepsilon}\|^2 - 1)$ are a variance/bias compromise

$$\frac{1}{n\sigma^2}\mathbb{E}[\|\mathbf{M}_R\underline{\varepsilon}\|^2] + \frac{1}{n\sigma^2}\|(\mathbf{I} - \mathbf{M}_R)\underline{r}_{\mathrm{add}}^{\mathrm{true}}\|^2,$$

which is minimized for $R = \infty$ and $h \propto n^{-1/5}$. Consequently for $\mathrm{AIC_T}$:

PROPOSITION 6. *Under the assumptions of Theorem 2 and $\mathbb{E}[\varepsilon^4] < \infty$, $\widehat{h}$ achieves the rate $n^{-1/5}$ and $\widehat{R} \to \infty$ (with probability tending to 1).*

If the true regression function is nonadditive, any $\widehat{R} \nrightarrow 0$ induces a bias of order $O(1)$, that is, an $\mathrm{AIC_T} - (\frac{1}{n\sigma^2}\|\underline{\varepsilon}\|^2 - 1)$ of $O(1)$. On the other hand, in well-behaved cases (continuous $\mathcal{S}_*^{-1}$), the optimal $\mathrm{AIC_T} - (\cdots)$ of the local linear estimator is of order $O(n^{-4/(4+d)})$, leading to $\widehat{R} \to 0$ and $\widehat{h} \propto n^{-1/(4+d)}$ in these cases.

4.6. *Investigating the rate of* $\widehat{R}$ *for* $\mathrm{AIC_T}$. Data-adaptive parameter selection is studied for fixed uniform design and additive regression functions; in this case, the penalty $\widehat{R}$ is large enough such that $\widehat{\underline{r}}_R - \widehat{\underline{r}}_{\mathrm{add}}$ becomes negligible.

As seen before, we can restrict ourselves to the case $h \propto n^{-1/5}$ and $R \to \infty$. In order to simplify the structure of $\widehat{r}_R^0(\underline{\mathbf{X}}_i)$, we assume a *fixed uniform design*: $\mathbf{S}(\underline{\mathbf{X}}_i)$ is diagonal and constant in the interior. Furthermore, we ignore boundary effects and $R$-dependency of $\mathcal{P}_{*,R}$. This allows us to simplify (13) as

$$(19) \qquad \widehat{r}_R^0(\underline{\mathbf{X}}_i) = \frac{1}{1+R}\widehat{r}_{ll}^0(\underline{\mathbf{X}}_i) + \frac{R}{1+R}\widehat{r}_{\mathrm{add}}^0(\underline{\mathbf{X}}_i).$$

By (19) we have $\mathbf{M}_R = \lambda\mathbf{M}_{ll} + (1-\lambda)\mathbf{M}_{\mathrm{add}}$ with $\lambda = \frac{1}{1+R}$. Hence, $\mathrm{AIC_T}$ is a polynomial of degree 2 in $\lambda$,

$$
\begin{aligned}
\mathrm{AIC_T} = &\ \lambda^2 \left\{ \frac{1}{n\sigma^2}\|\mathbf{M}_{ll}\,\underline{\varepsilon}\|^2 + \frac{1}{n\sigma^2}\|\mathbf{M}_{\mathrm{add}}\,\underline{\varepsilon}\|^2 - \frac{2}{n\sigma^2}\langle\mathbf{M}_{ll}\,\underline{\varepsilon}, \mathbf{M}_{\mathrm{add}}\,\underline{\varepsilon}\rangle \right\} \\
&+ \lambda \left\{ \frac{2}{n\sigma^2}\mathbb{E}[\langle\mathbf{M}_{ll}\,\underline{\varepsilon}, \mathbf{M}_{\mathrm{add}}\,\underline{\varepsilon}\rangle] - \frac{2}{n\sigma^2}\|\mathbf{M}_{\mathrm{add}}\,\underline{\varepsilon}\|^2 \right. \\
&\qquad\quad \left. + \frac{2}{n\sigma^2}(1-\mathbb{E})[\langle\underline{\varepsilon}, \mathbf{M}_{ll}\,\underline{\varepsilon}\rangle + \langle\mathbf{M}_{ll}\,\underline{\varepsilon}, \mathbf{M}_{\mathrm{add}}\,\underline{\varepsilon}\rangle - \langle\underline{\varepsilon}, \mathbf{M}_{\mathrm{add}}\,\underline{\varepsilon}\rangle] \right\} \\
&+ 1 \left\{ \frac{1}{n\sigma^2}\|(\mathbf{I}-\mathbf{M}_R)\,\underline{r}_{\mathrm{add}}^{\mathrm{true}}\|^2 + O_P(h^2/\sqrt{n}) + \left(\frac{1}{n\sigma^2}\|\underline{\varepsilon}\|^2 - 1\right) \right. \\
&\qquad\quad \left. + \frac{1}{n\sigma^2}\|\mathbf{M}_{\mathrm{add}}\,\underline{\varepsilon}\|^2 + \frac{2}{n\sigma^2}(1-\mathbb{E})[\langle\underline{\varepsilon}, \mathbf{M}_{\mathrm{add}}\,\underline{\varepsilon}\rangle] \right\}.
\end{aligned}
$$

In Appendix A.0.12 we show that the dominating terms of $\lambda^2$ and $\lambda$ are $\frac{1}{n\sigma^2}\mathbb{E}[\|\mathbf{M}_{ll}\,\underline{\varepsilon}\|^2] \propto \frac{1}{nh^d}$ and $(1-\mathbb{E})\langle\underline{\varepsilon}, \mathbf{M}_{ll}\,\underline{\varepsilon}\rangle = O_P(\frac{1}{nh^{d/2}})$, respectively, leading to $\lambda_{\min} \approx -(1-\mathbb{E})\langle\underline{\varepsilon}, \mathbf{M}_{ll}\,\underline{\varepsilon}\rangle / \mathbb{E}[2\|\mathbf{M}_{ll}\,\underline{\varepsilon}\|^2] = O_P(h^{d/2})$.

PROPOSITION 7. *Under the assumptions of Theorem 2, fixed uniform design and* $\mathbb{E}[\varepsilon^4] < \infty$, *we obtain* $\widehat{R}^{-1} = O_P(n^{-d/10})$ *for* $\mathrm{AIC_T}$.

(See Appendix A.0.12 for a proof.) By Theorem 2, $\widehat{R}$ grows fast enough to ensure $\|\widehat{\underline{r}}_R - \widehat{\underline{r}}_{\mathrm{add}}\|_2^2 = o_p(n^{-4/5})$.

Note that the rate of the minimum is not affected by the $\mathbf{M}_{\mathrm{add}}$-dependent terms. Accordingly, assuming $\mathcal{P}_{*,R} = \mathcal{P}_*$ is not critical.

*Comparison of* $\mathrm{AIC_C}$ *and* $\mathrm{AIC_T}$. By (19), $\mathrm{tr}(\mathbf{M}_R) = \lambda\,\mathrm{tr}(\mathbf{M}_{ll}) + (1-\lambda) \times \mathrm{tr}(\mathbf{M}_{\mathrm{add}})$ is monotone decreasing in $R$, because (asymptotically) $1 > \frac{1}{n}\,\mathrm{tr}(\mathbf{M}_{ll}) > \frac{1}{n}\,\mathrm{tr}(\mathbf{M}_{\mathrm{add}})$. If we add $(\frac{1}{n}\,\mathrm{tr}(\mathbf{M}_R))^{2+\ell}$ ($\ell \geq 0$) to $\mathrm{AIC_T}$, $\lambda_{\min}$ becomes smaller. Hence, $\mathrm{AIC_C}$ [with $\log(\widehat{\sigma}^2)$ replaced by $\log(\sigma^2) + \frac{\widehat{\sigma}^2}{\sigma^2} - 1$] chooses $\widehat{R}$ at least as large as $\text{AIC}_\text{T}$. However, this effect is asymptotically negligible as the leading terms are unchanged.

**5. Finite sample evaluation.** For the example in Section 1, we compared the penalized estimator with the *local linear* ($R = 0$) and the *additive* ($R = \infty$) estimator and obtained a lower integrated squared error (ISE) for the penalized estimator. As seen in Figure 3, data-adaptive choice (specified in Section 4.5) of the parameters $R$ and $h$ is successful: the theoretical improvement due to generalization holds also in practice. In Section 5.1, we will see whether this holds for other situations. Furthermore, we will investigate how the estimator performs in the special case of an additive model; see Section 5.2. Finally (Section 5.3), we apply our estimator to the ozone dataset already analyzed by Hastie and Tibshirani (1990).

5.1. *Nonadditive regression function.* We will examine 50 realizations of the same kind as in Section 1. Later on, we summarize the effect of a nonuniform design density and a larger sample size.

In the example in Section 1, the optimal penalty parameter is larger than zero and the penalized estimator outperforms the local linear estimator when using optimal parameters. The optimal parameters are approximated sufficiently well by $\text{AIC}_\text{C}$.

Here we generate 50 realizations of the data as follows: the design consists of 200 random observations $\underline{\mathbf{X}}_i$, uniform in $[0, 1]^2$. The response $Y_i$ is $r^{\text{true}}(\underline{\mathbf{X}}_i) + \varepsilon_i$, where $r^{\text{true}}$ is shown in Figure 1(b) [$\underline{\mathbf{1}} = (1, 1)^\top$]:

$$(20) \quad r^{\text{true}}(\underline{\mathbf{x}}) = 15e^{-32\|\underline{\mathbf{x}} - (1/4)\underline{\mathbf{1}}\|^2} + 35e^{-128\|\underline{\mathbf{x}} - (3/4)\underline{\mathbf{1}}\|^2} + 25e^{-2\|\underline{\mathbf{x}} - (1/2)\underline{\mathbf{1}}\|^2}$$

and $\varepsilon_i$ is normally distributed ($\sigma = 5$).

In order to find the optimal parameters for each realization of the design, we calculated the ISE for different pairs $(R, h)$ [see Figure 3(a)] and performed a grid search. Actually, $(\frac{R}{1+R}, \log_{10}(h))$ is equidistant with resolution $(0.01, 0.005)$. Similarly, we find the minimizers of $\text{AIC}_\text{C}$ and GCV.

The simulation is summarized by 50 realizations of ISE evaluated for the *penalized*, the *local linear* and the *additive* estimator using *optimal* and *data-adaptive* parameters.

Let us introduce some notation. The global minimum at $(R_{\text{opt}}, h_{\text{opt}})$ is ISE(opt). The minimum of the local linear estimator ($R = \frac{1}{9999}$ instead of 0 for numerical reasons) is ISE(opt, $R = 0$) and the minimum of the additive estimator ($R = 9999$ instead of $\infty$) is ISE(opt, $R = \infty$). Data-adaptive parameters $(R_{\text{AIC}}, h_{\text{AIC}})$ are obtained by finding the minimizer of $\text{AIC}_\text{C}$ (Section 4.5) on a grid. The corresponding ISE is denoted by ISE(AIC). Analogously, we write ISE(AIC, $R = 0$) for the local linear and ISE(AIC, $R = \infty$) for the additive estimator.

Table 2 and Figure 4 summarize the results of this evaluation. Given the optimal parameter values for $R$ and $h$, penalized estimation has clearly the potential for improvement compared to fitting the full model with a median percentage gain of 17% [item (a) of Table 2 and Figure 4]. This relative gain is larger for realizations with a small ISE. The additive estimator is not competitive and will hence not be shown.

To be able to achieve these gains in practice, we need a good method for parameter selection. The corrected Akaike criterion $AIC_C$ is such a method, and is moreover computationally attractive. When comparing (c) of Table 2 with (a), we see that the performance based on estimated $R$ and $h$ is almost as good as that based on optimal parameters. Item (b) shows that data-adaptive parameter selection via $AIC_C$ is attractive: a median increase in relative ISE of only 10% has to be tolerated.

Interestingly, application of the full model with optimal bandwidths is clearly inferior to using the penalized estimator with data-driven parameter selection [see item (d)].

*Other situations.* The above simulation was also carried out for two *nonuniform designs* on $[0, 1]^2$,

$$f_1(x_1, x_2) = \tfrac{1}{2} + \tfrac{1}{2}(x_1 + x_2) \quad \text{and} \quad f_2(x_1, x_2) = \tfrac{3}{2} - \tfrac{1}{2}(x_1 + x_2).$$

Both are linear in $(x_1 + x_2)$ and have range $(0.5, 1.5)$. Density $f_1$ is preferred over $f_2$ because of the high peak in $r^{\text{true}}$ at $(0.75, 0.75)$. This is reflected by the optimal penalty $R_{\text{opt}}$: compared with the uniform design, $f_1$ needs a larger and $f_2$ a smaller penalty; see Table 3(e). Similarly, the ideal relative gain due to penalizing is larger for $f_1$ and smaller for $f_2$, item (a). Again, the performance remains almost as good when selecting parameters $R$ and $h$ via $AIC_C$ [see item (c)]. The penalized estimator with $AIC_C$-selected parameters is better than the local linear estimator with optimal parameters; however, for density $f_2$ the difference becomes smaller.

TABLE 2

*Quantiles for ideal relative gain due to penalizing* (a), *for loss due to* $AIC_C$ *selection* (b), *for relative gain (loss) due to penalizing for data-adaptive parameters* (c), *for relative difference between data-adaptive penalized and optimal full estimator* (d) *and for $R_{\text{opt}}$* (e)

| Model defined in (20) | | min | 10% | med | 90% | max |
|---|---|---|---|---|---|---|
| (a) | $\frac{\text{ISE(opt}, R=0)-\text{ISE(opt)}}{\text{ISE(opt)}}$ | 1.8% | 4.9% | 17% | 49% | 65% |
| (b) | $\frac{\text{ISE(AIC)}-\text{ISE(opt)}}{\text{ISE(opt)}}$ | 0.2% | 1.1% | 10% | 24% | 66% |
| (c) | $\frac{\text{ISE(AIC}, R=0)-\text{ISE(AIC)}}{\text{ISE(AIC)}}$ | −8.0% | 1.3% | 16% | 47% | 149% |
| (d) | $\frac{\text{ISE(opt}, R=0)-\text{ISE(AIC)}}{\text{ISE(opt)}}$ | −50% | −9.7% | 11% | 27% | 51% |
| (e) | $R_{\text{opt}}$ | 0.03 | 0.06 | 0.12 | 0.25 | 0.32 |

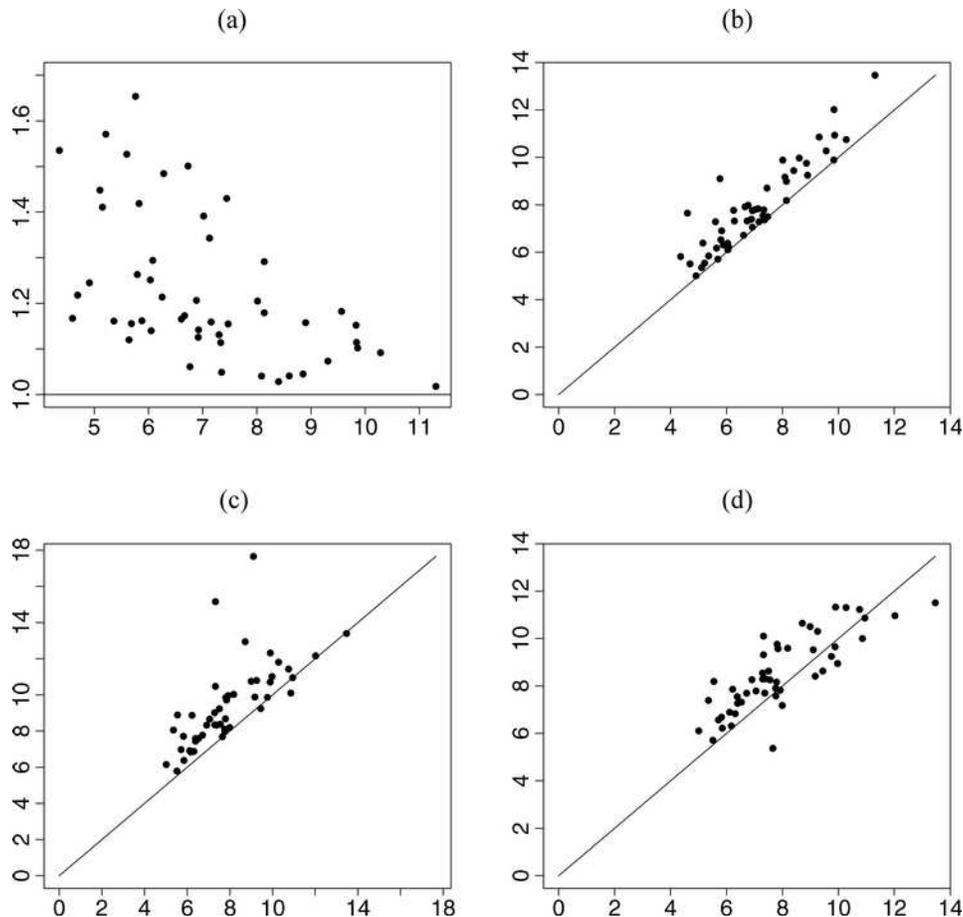

Fig. 4. *Comparison of ISE performance for the nonadditive regression function* (20). *Relative gain due to penalizing depending on* ISE(opt) (a) $[y = \text{ISE(opt}, R = 0)/\text{ISE(opt)}$ *vs.* $x = \text{ISE(opt)}]$, *effect of* $\text{AIC}_C$ *selection* (b) $[y = \text{ISE(AIC)}$ *vs.* $x = \text{ISE(opt)}]$, *comparison of penalized vs. full modeling, parameters data-driven* (c) $[y = \text{ISE(AIC}, R = 0)$ *vs.* $x = \text{ISE(AIC)}]$, *comparison of penalized modeling ( parameters data-adaptive) with full modeling (optimal bandwidth)* (d) $[y = \text{ISE(opt}, R = 0)$ *vs.* $x = \text{ISE(AIC)}]$.

When doubling $n$ to 400, parameter selection is improved; see Table 3, column "400," item (b). Because of the smaller $R_{\text{opt}}$, the gain due to penalizing is smaller but still not negligible.

Parameter selection by GCV instead of $\text{AIC}_C$ shows the same pattern (data not shown).

5.2. *Additive regression function.* For additive regression functions, the question arises whether we pay a price for the additional flexibility of penal-

TABLE 3
*For different design densities we compare the medians of the same quantities as in Table 2*

| Model defined in (20) | | $f_1$ | unif | $f_2$ | 400 |
|---|---|---|---|---|---|
| (a) | $\frac{\text{ISE(opt}, R=0) - \text{ISE(opt)}}{\text{ISE(opt)}}$ | 19% | 17% | 6% | 12% |
| (b) | $\frac{\text{ISE(AIC}) - \text{ISE(opt)}}{\text{ISE(opt)}}$ | 10% | 10% | 4.3% | 4.5% |
| (c) | $\frac{\text{ISE(AIC}, R=0) - \text{ISE(AIC)}}{\text{ISE(opt)}}$ | 20% | 16% | 11% | 13% |
| (d) | $\frac{\text{ISE(opt}, R=0) - \text{ISE(AIC)}}{\text{ISE(opt)}}$ | 11% | 11% | 1.5% | 6.7% |
| (e) | $R_{\text{opt}}$ | 0.15 | 0.12 | 0.10 | 0.09 |

The number of realizations is 50, the sample size is $n = 200$ and the name of the column denotes the design density—except the last column ($n = 400$, uniform).

izing local linear estimation compared with additive estimation. Therefore, we choose an additive regression function and examine 50 realizations.

We generated data using the regression function

$$r^{\text{true}}(\underline{\mathbf{x}}) = \sum_{k=1}^{2} \left( \tfrac{15}{2} e^{-32(x_k - 1/4)^2} + \tfrac{35}{2} e^{-128(x_k - 3/4)^2} + \tfrac{25}{2} e^{-2(x_k - 1/2)^2} \right).$$

Uniform design $\underline{\mathbf{X}}_i$ and errors $\varepsilon_i$ ($\sigma = 5$) are the same as in Section 5.1. Estimating the additive model can be considered easy, as the data are rich enough for multivariate local linear estimation.

Since $\text{AIC}_C$ has no problems with undersmoothing, we ignore in the simulations the impracticable condition in Section 4.5—excluding undersmoothing—which was imposed for classical AIC. For $\text{AIC}_C$ the additive model is detected in 47 out of 50 cases, as $R_{\text{AIC}}$ attains the maximal value. In the remaining three realizations we obtained a relative loss in ISE of 0.6% and 5.2% in two cases; a gain of 3.5% was achieved in one case. Hence, *model choice* by $\text{AIC}_C$ was successful in this example.

Model selection by GCV detected the additive model in 24/50 cases only, whereas classical AIC failed completely (0/50).

5.3. *Application to ozone data.* We apply our method to the ozone dataset using three out of nine predictors. The penalized estimator detects relevant deviations from an additive model. The local linear estimator produces artifacts, which do not occur in the penalized estimator.

We used the ozone dataset from the R package gss; see Hastie and Tibshirani [(1990), Section 10.3]. The variable "wind speed" (wdsp) contains one excessive value (observation number 92) which was removed, leading to $n = 329$. The dependent variable $\underline{\mathbf{Y}}$ was chosen as the logarithm of the "upland ozone concentration" (upo3). Using gam (package mgcv), we chose those three predictors which maximize adjusted R-squared among additive models with *bivariate interaction terms* with 16 degrees of freedom each: "humidity" (hmdt), "inversion base height" (ibtp), and "calendar day" (day).

Note that this additive model with *bivariate interaction terms* has roughly the same adjusted R-squared as the additive model with all nine predictors and four degrees of freedom for each component; see Table 4. Hence, when using these three predictors, we expect substantial information in the interaction terms.

The three variables in this model were scaled to $[0, 1]$. Let univariate bandwidths $h_1$, $h_2$ and $h_3$ correspond to four degrees of freedom each, as in Hastie and Tibshirani (1990). These bandwidths lie close together (min = 95% max) with mean $h = \sqrt[3]{h_1 h_2 h_3}$ at 0.237. Parameters $R$ and $c$ are selected by AIC$_C$, such that the bandwidths are $ch_1$, $ch_2$ and $ch_3$, respectively.

For the penalized estimator, AIC$_C$ selected $R = 0.04$ and $ch = 0.2065$ and is clearly nonadditive. For the local linear estimator, AIC$_C$ selected $ch = 0.240$. The lower half of Table 4 demonstrates that the penalized estimator is vastly better than the additive one in terms of adjusted $R$-squared, and slightly better than the local linear estimator.

Next, we orthogonally decompose the local linear and the penalized estimator into intercept, additive components, bivariate interactions and remainder. Penalizing shrinks the bivariate interactions and the remainder; see Table 5.

Figure 5 compares the local linear and penalized estimators, univariate components on the top and the largest bivariate interaction (ibtp, hmdt) on the bottom. The plots for univariate components demonstrate those regularization properties of the penalized estimator. The plot in the center of the bottom row shows the design and the smoothing windows for one

TABLE 4
*Adjusted R-squared for different models and estimators*

| Estimator | | # independent variables | $R$-squared (adj.) |
|---|---|---|---|
| Regression | additive | 9 | 82.5% |
| spline (gam) | additive | 3 | 73.9% |
| $df = 4 \mid 16$ | additive + bivariate interaction | 3 | 81.3% |
| Penalized | $\widehat{r}_{ll}$ multivariate ($R = 10^{-4}$) | 3 | 81.7% |
| local linear | $\widehat{r}_R$ penalized ($R = 0.04$) | 3 | 82.5% |
| AIC$_C$ | $\widehat{r}_{\text{add}}$ additive ($R = \infty$) | 3 | 73.9% |

Above, we use regression splines with a fixed number of knots. Below, we use local linear estimators with AIC$_C$ selected parameters. The two additive estimators with three predictors are equivalent to each other but inferior to the rest.

*Orthogonal decomposition of estimation on a grid into constant, additive, interaction and remainder components*

| $R$ | $h$ | $r_0$ | $r_1$ | $r_2$ | $r_3$ | $r_{12}$ | $r_{13}$ | $r_{23}$ | $r_{123}$ |
|---|---|---|---|---|---|---|---|---|---|
| $10^{-4}$ | 0.240 | 3.95 | 0.482 | 0.054 | 0.025 | 0.023 | 0.082 | 0.019 | 0.059 |
| 0.04 | 0.207 | 3.96 | 0.478 | 0.051 | 0.017 | 0.009 | 0.023 | 0.002 | 0.011 |

We compare the mean squares of each component for the penalized ($R = 0.04$) and for the local linear ($R = 10^{-4}$) estimator. Penalizing shrinks interaction and remainder components. Univariate additive components are slightly reduced. The indices are $1 = \text{ibtp}$, $2 = \text{day}$ and $3 = \text{hmdt}$.

output point. Keep in mind that the smoothing windows are actually three-dimensional cubes and hence not all points inside the rectangle actually contribute to the local linear estimator.

Let us mention that parameter selection criteria such as $\text{AIC}_C$ and GCV evaluate the estimator at design points and hence are not influenced by its behavior in sparse regions. Comparing the local linear estimator with the penalized estimator, we observe some strange structure for the former at hmdt = 0.3. This is clearly an artifact, as for day $\geq 0.75$ and ibtp $\geq 0.75$, the local linear estimator is an *extrapolation*.

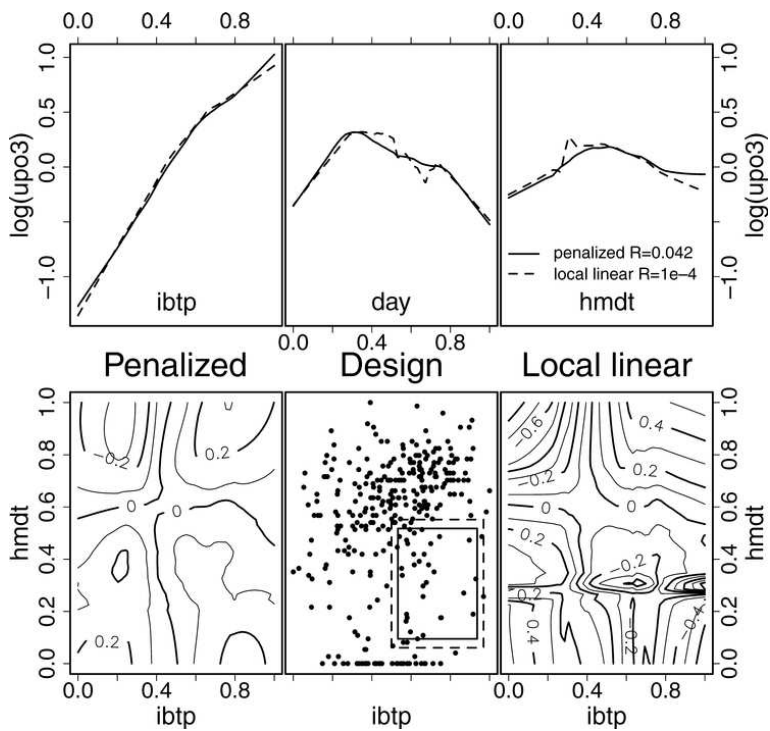

FIG. 5. *Comparison of penalized and local linear estimators. Above, the univariate additive components are shown. Below, the bivariate components of ibtp and hmdt are compared. In between, the design is shown including the smoothing windows at $(0.73, 0.31)$.*

We conclude that the penalized estimator outperforms the additive estimator and is also superior to the full estimator regarding adjusted $R$-squared and regularization properties.

*Reproducing simulation results.* An implementation of $\widehat{\underline{r}}_R$ together with the R code used in the simulations of this paper is provided at www.biostat.unizh.ch/Software.

## APPENDIX

### A.0.1. Assumptions and details.

A.0.2. *Conditions for optimality of the MLN estimator.* MLN show that the estimator $\widehat{\underline{r}}_{\text{add}}$ is asymptotically equal to the oracle estimator if $(Y_i, \underline{\mathbf{X}}_i)$ are i.i.d., the true regression function $r^{\text{true}}(\underline{\mathbf{X}}_i) = \mathbb{E}[Y_i|\underline{\mathbf{X}}_i]$ is additive (1) and the following conditions hold:

CONDITION MLN:B1. The kernel $K$ is bounded, has compact support ($[-C_1, C_1]$), is symmetric about zero and is Lipschitz continuous.

CONDITION MLN:B2$'$. The $d$-dimensional vector $\underline{\mathbf{X}}$ has compact support $[0,1]^d$ and its density $f$ is bounded away from zero and infinity on $[0,1]^d$.

The product kernel $K_h$ with bandwidths $h_1, \ldots, h_d$ is constructed from the univariate kernel $K$ by $K_h(\underline{\mathbf{X}}, \underline{\mathbf{x}}) = \prod_{k=1}^d K([\underline{\mathbf{X}} - \underline{\mathbf{x}}]_k/h_k)/h_k$.

Furthermore, the kernel is rescaled at the boundary such that for all $\underline{\mathbf{X}}_i \in [0,1]^d$,

$$\int_{[0,1]^d} K_h(\underline{\mathbf{X}}_i, \underline{\mathbf{x}}) \, d\underline{\mathbf{x}} = 1.$$

This modification does not affect the local linear estimator, but it changes its *projection* to the additive model. Hence, the estimation of the marginal design density is equal to an integrated full-dimensional density estimation. Additionally to MLN, we assume $K \geq 0$.

CONDITION MLN:B3$'$. For some $\theta > 5/2$, $\mathbb{E}[|Y|^\theta] < \infty$.

Additionally to MLN, we assume $\mathbb{E}[\varepsilon^4] < \infty$ in Sections 4.5 and 4.6.

CONDITION MLN:B4′.   The true regression function $r^{\text{true}}(\underline{\mathbf{x}}) = \mathbb{E}\,[Y \,|\, \underline{\mathbf{X}} = \underline{\mathbf{x}}]$ is twice continuously differentiable and $f$ is once continuously differentiable.

CONDITION C1.   Assume there exist constants $c_k$ with $n^{1/5}h_k \to c_k$, $k = 1, \dots, d$.

CONDITION C1+.   The bandwidths $h_1, \dots, h_d$ are as in Condition C1 or larger. As a matter of course we assume that $h_k \to 0$.

A.0.3. *Definition of* $\mathbf{Z}$.   In Section 2.2 the output grid $\underline{\mathbf{t}}_j$, $j = 1, \dots, m$, and the parameters $\underline{\boldsymbol{\beta}} = \text{col}_j(\underline{\boldsymbol{\beta}}^{(j)}) = \text{col}_j(\text{col}_\ell(\beta_\ell^{(j)}))$ were introduced. In Section 3.2 $\mathbf{P}_{\text{add}}$ was decomposed into the product $\mathbf{Z}^\top\mathbf{Z}$. Instead of writing down the matrix $\mathbf{Z}$, we show what $\mathbf{Z}$ does with a vector $\underline{\boldsymbol{\beta}} \in \mathbf{F}_{\text{full}}$. For the index ranges we use $j = 1, \dots, m$; $\ell = 0, \dots, d$; $k = 1, \dots, d$. Let $\underline{\boldsymbol{\beta}}_\ell = \text{col}_j(\beta_\ell^{(j)})$. Define

$$\mathbf{Z}\underline{\boldsymbol{\beta}} = \text{col}(\mathbf{Z}_{01}\underline{\boldsymbol{\beta}}_0, \mathbf{Z}_2\underline{\boldsymbol{\beta}}_0, \dots, \mathbf{Z}_d\underline{\boldsymbol{\beta}}_0, \mathbf{Z}_{01}\underline{\boldsymbol{\beta}}_1, \mathbf{Z}_{02}\underline{\boldsymbol{\beta}}_2, \dots, \mathbf{Z}_{0d}\underline{\boldsymbol{\beta}}_d),$$

where $\mathbf{Z}_{0k}\underline{\boldsymbol{\beta}}_k$ adds those $\beta_k^{(j)}$ which have the $k$th coordinate of $\underline{\mathbf{t}}_j$ in common,

$$\mathbf{Z}_{0k}\underline{\boldsymbol{\beta}}_k = \sqrt{\frac{m_k}{m}} \sum_j \beta_k^{(j)} \begin{pmatrix} 1_{t_1^k = t_{j,k}} \\ \vdots \\ 1_{t_{m_k}^k = t_{j,k}} \end{pmatrix}.$$

For identifiability, all additive components of the intercept except the first one should have mean zero. Therefore, $\mathbf{Z}_k\underline{\boldsymbol{\beta}}_0$ is defined as $\mathbf{Z}_{0k}\underline{\boldsymbol{\beta}}_0$ with subtracted mean.

We did not modify $\mathbf{Z}_k$ to have full rank, as this makes implementation more complicated, and it appears that the additional computing steps offset the gain due to lower dimension. $\mathbf{Z}^\top\mathbf{Z}$ is a projection and hence $\mathbf{Z}\mathbf{Z}^\top$ is too.

### A.0.4. Proofs.

A.0.5. *Algorithm, structure and proofs for Section* 3.2.

*Deriving* (7).  Rao and Kleffe (1988) provide a *generalized inverse* for the matrix $\mathbf{B} + \mathbf{C}\mathbf{D}\mathbf{C}^\top$,

$$\mathbf{B}^- - \mathbf{B}^-\mathbf{C}\mathbf{D}(\mathbf{I} + \mathbf{C}^\top\mathbf{B}^-\mathbf{C}\mathbf{D})^-\mathbf{C}^\top\mathbf{B}^-.$$

This holds if $\mathbf{B}$ and $\mathbf{D}$ are symmetric and if the image of $\mathbf{B}$ contains the image of $\mathbf{C}$. We apply this formula for $(\mathbf{B}, \mathbf{C}, \mathbf{D}) = (\mathbf{S} + R\mathbf{I}, \mathbf{Z}^\top, -R\mathbf{I})$.

*Verify that* $\mathbf{P}_{\mathbf{S},R}$ *is a projection from* $\mathbf{F}_{\text{full}}$ *to* $\mathbf{F}_{\text{add}}$. Define

$$(21) \qquad \Lambda_{\mathbf{S},R} := \mathbf{I} - \mathbf{Z}\mathbf{A}_R\mathbf{Z}^\top = (\mathbf{I} - \mathbf{Z}\mathbf{Z}^\top) + \mathbf{Z}(\mathbf{I} - \mathbf{A}_R)\mathbf{Z}^\top.$$

Because $\mathbf{I} - \mathbf{A}_R \geq \mathbf{0}$, the image of $\Lambda_{\mathbf{S},R}$ contains the image of $\mathbf{Z}(\mathbf{I} - \mathbf{A}_R)$. We get

$$(22) \qquad \Lambda_{\mathbf{S},R}\Lambda_{\mathbf{S},R}^-\mathbf{Z}(\mathbf{I} - \mathbf{A}_R) = \mathbf{Z}(\mathbf{I} - \mathbf{A}_R).$$

Furthermore,

$$(23) \qquad (\mathbf{I} - \mathbf{Z}\mathbf{Z}^\top)\mathbf{Z} = \mathbf{0}$$

and the definition of $\Lambda_{\mathbf{S},R}$ implies

$$(24) \qquad (\mathbf{I} - \mathbf{Z}\mathbf{Z}^\top)\Lambda_{\mathbf{S},R} = (\mathbf{I} - \mathbf{Z}\mathbf{Z}^\top).$$

Applying (22)–(24) leads to

$$(25) \qquad (\mathbf{I} - \mathbf{Z}\mathbf{Z}^\top)\Lambda_{\mathbf{S},R}^-\mathbf{Z}(\mathbf{I} - \mathbf{A}_R) = \mathbf{0}.$$

$\mathbf{P}_{\mathbf{S},R} = \mathbf{Z}^\top\Lambda_{\mathbf{S},R}^-\mathbf{Z}(\mathbf{I} - \mathbf{A}_R)$ is a *projection* because

$$
\begin{aligned}
\mathbf{P}_{\mathbf{S},R}^2 &= \mathbf{Z}^\top\Lambda_{\mathbf{S},R}^-\mathbf{Z}(\mathbf{I} - \mathbf{A}_R)\mathbf{Z}^\top\Lambda_{\mathbf{S},R}^-\mathbf{Z}(\mathbf{I} - \mathbf{A}_R)\\[4pt]
&\stackrel{(21)}{=} \mathbf{Z}^\top\Lambda_{\mathbf{S},R}^-(\Lambda_{\mathbf{S},R} - (\mathbf{I} - \mathbf{Z}\mathbf{Z}^\top))\Lambda_{\mathbf{S},R}^-\mathbf{Z}(\mathbf{I} - \mathbf{A}_R)\\[4pt]
&\stackrel{(25)}{=} \mathbf{Z}^\top\Lambda_{\mathbf{S},R}^-\Lambda_{\mathbf{S},R}\Lambda_{\mathbf{S},R}^-\mathbf{Z}(\mathbf{I} - \mathbf{A}_R)\\[4pt]
&\stackrel{(22)}{=} \mathbf{Z}^\top\Lambda_{\mathbf{S},R}^-\mathbf{Z}(\mathbf{I} - \mathbf{A}_R) = \mathbf{P}_{\mathbf{S},R}
\end{aligned}
$$

and $(\mathbf{I} - \mathbf{A}_R)\mathbf{P}_{\mathbf{S},R}$ is symmetric. If $\Lambda_{\mathbf{S},R}$ is nonsingular, the image of $\mathbf{P}_{\mathbf{S},R}$ is $\mathbf{F}_{\text{add}}$. Let us verify that $\mathbf{P}_{\mathbf{S},R}\mathbf{P}_{\text{add}} = \mathbf{P}_{\text{add}}$:

$$
\begin{aligned}
(\mathbf{I} - \mathbf{P}_{\mathbf{S},R})\mathbf{Z}^\top\mathbf{Z} &= \mathbf{Z}^\top(\mathbf{I} - \Lambda_{\mathbf{S},R}^{-1}\mathbf{Z}(\mathbf{I} - \mathbf{A}_R)\mathbf{Z}^\top)\mathbf{Z}\\[4pt]
&\stackrel{(21)}{=} \mathbf{Z}^\top(\mathbf{I} - \Lambda_{\mathbf{S},R}^{-1}(\Lambda_{\mathbf{S},R} - (\mathbf{I} - \mathbf{Z}\mathbf{Z}^\top)))\mathbf{Z}\\[4pt]
&\stackrel{(23)}{=} \mathbf{Z}^\top(\mathbf{I} - \Lambda_{\mathbf{S},R}^{-1}\Lambda_{\mathbf{S},R})\mathbf{Z} = \mathbf{0}.
\end{aligned}
$$

Formula (8) is straightforward.

A.0.6. *Iterative formula and proofs for Section* 3.3.

*Convergence of iterative algorithm in* (9) *in the case of nonuniqueness.* We denote by $(\mathbf{Z}\mathbf{A}_R\mathbf{Z}^\top)^\infty$ the projection to the subspace of eigenvectors of $\mathbf{Z}\mathbf{A}_R\mathbf{Z}^\top$ with eigenvalues 1. Because $\mathbf{Z}(\mathbf{I} - \mathbf{A}_R)\widehat{\underline{\beta}}_{ll}$ is orthogonal to the above subspace, $(\mathbf{Z}\mathbf{A}_R\mathbf{Z}^\top)^\infty\underline{\gamma}^{[a]} = (\mathbf{Z}\mathbf{A}_R\mathbf{Z}^\top)^\infty\underline{\gamma}^{[0]}$ and $\mathbf{I} - (\mathbf{Z}\mathbf{A}_R\mathbf{Z}^\top)^\infty\underline{\gamma}^{[a]}$ converges.

*Iteration formula for large R, deriving* (10). Starting with (8), multiply by $\mathbf{Z}$, replace $\mathbf{P}_{\mathbf{S},R}\underline{\hat{\boldsymbol{\beta}}}_{ll}$ by $\mathbf{Z}^{\top}\mathbf{Z}\underline{\hat{\boldsymbol{\beta}}}_R$, subtract $\mathbf{Z}\underline{\hat{\boldsymbol{\beta}}}_R$, multiply by $\alpha R$, add $\mathbf{Z}\underline{\hat{\boldsymbol{\beta}}}_R$, and apply $(\mathbf{I} - \mathbf{Z}\mathbf{Z}^{\top})\mathbf{Z}\underline{\hat{\boldsymbol{\beta}}}_R = \underline{\mathbf{0}}$ because of (23).

A.0.7. *Properties on grid and proofs for Section* 4.1. Proposition 1 is an interpretation of (8) and requires no proof.

PROOF OF PROPOSITION 2. We need to prove that $\mathbf{P}_{\mathrm{add}}\underline{\hat{\boldsymbol{\beta}}}_R = \mathbf{P}_{\mathbf{S},R}\underline{\hat{\boldsymbol{\beta}}}_{ll}$, respectively

$$\mathbf{P}_{\mathrm{add}}((\mathbf{I} - \mathbf{A}_R)(\mathbf{I} - \mathbf{P}_{\mathbf{S},R})) = \mathbf{0}.$$

Transposing and applying the symmetry of the matrix, this is equivalent to

$$(\mathbf{I} - \mathbf{A}_R)(\mathbf{I} - \mathbf{P}_{\mathbf{S},R})\mathbf{P}_{\mathrm{add}} = \mathbf{0}.$$

This holds because $\mathbf{P}_{\mathbf{S},R}\mathbf{P}_{\mathrm{add}} = \mathbf{P}_{\mathrm{add}}$ [if and only if (5) is unique]. $\square$

A.0.8. *Definition of* $\mathcal{P}_{*,R}$ *and proofs for Section* 4.2. We use the abbreviation $\mathcal{S}_{\mathrm{add}} = \mathcal{P}_{\mathrm{add}}\mathcal{S}_*\mathcal{P}_{\mathrm{add}}$ restricted to $\mathcal{F}_{\mathrm{add}}$. Hence, $\mathcal{S}_{\mathrm{add}}$ is a linear operator on $\mathcal{F}_{\mathrm{add}}$, and it has a continuous inverse with probability tending to 1 for $n \to \infty$; see Lemma 2 in Section 4.4.

Define the operator $\Lambda_{*,R} : \mathcal{F}_{\mathrm{add}} \to \mathcal{F}_{\mathrm{add}}$:

$$\Lambda_{*,R} = \mathcal{P}_{\mathrm{add}}(\mathcal{S}_* + R\mathcal{I})^{-1}\mathcal{S}_*\mathcal{P}_{\mathrm{add}}.$$

If $\mathcal{S}_{\mathrm{add}}^{-1}$ exists and is continuous, $\Lambda_{*,R}^{-1}$ is continuous because

$$\Lambda_{*,R} \geq \frac{1}{\|\mathcal{S}_*\|_{2,\sup} + R}\mathcal{S}_{\mathrm{add}}.$$

Let us define the projection $\mathcal{P}_{*,R} : \mathcal{F}_{\mathrm{full}} \to \mathcal{F}_{\mathrm{add}}$ by

$$(26) \qquad \mathcal{P}_{*,R} = \mathcal{P}_{\mathrm{add}}\Lambda_{*,R}^{-1}\mathcal{P}_{\mathrm{add}}(\mathcal{S}_* + R\mathcal{I})^{-1}\mathcal{S}_*.$$

$\mathcal{S}_*$ is continuous ($\|\mathcal{S}_*\|_{2,\sup} < \infty$), because kernel weights are bounded and have compact support (Condition MLN:B1).

Below we will verify that the choice for $\underline{\hat{r}}_R$ in (13) in Section 4.2 satisfies the normal equation (12). We need the properties $\mathcal{P}_{\mathrm{add}}\mathcal{P}_{*,R} = \mathcal{P}_{*,R}$ and

$$
\begin{aligned}
(27) \qquad \mathcal{P}_{\mathrm{add}}(\mathcal{I} - (\mathcal{S}_* + R\mathcal{I})^{-1}R)\mathcal{P}_{*,R} &\overset{(26)}{=} \Lambda_{*,R}\Lambda_{*,R}^{-1}\mathcal{P}_{\mathrm{add}}(\mathcal{S}_* + R\mathcal{I})^{-1}\mathcal{S}_* \\
&= \mathcal{P}_{\mathrm{add}}(\mathcal{S}_* + R\mathcal{I})^{-1}\mathcal{S}_*.
\end{aligned}
$$

We have to verify the normal equations $(\mathcal{S}_* + R(\mathcal{I} - \mathcal{P}_{\mathrm{add}}))\widehat{\underline{r}}_R = \mathcal{S}_* \widehat{\underline{r}}_{ll}$:

$$
\begin{aligned}
(\mathcal{S}_* \; + \; & R(\mathcal{I} - \mathcal{P}_{\mathrm{add}}))\widehat{\underline{r}}_R \\
&= ((\mathcal{S}_* + R\mathcal{I}) - R\mathcal{P}_{\mathrm{add}})(\mathcal{S}_* + R\mathcal{I})^{-1}\{\mathcal{S}_* + R\mathcal{P}_{*,R}\}\widehat{\underline{r}}_{ll} \\
&= \; \mathcal{S}_* \widehat{\underline{r}}_{ll} + R\mathcal{P}_{\mathrm{add}}\{(I - (\mathcal{S}_* + R\mathcal{I})^{-1}R)\mathcal{P}_{*,R} - (\mathcal{S}_* + R\mathcal{I})^{-1}\mathcal{S}_*\}\widehat{\underline{r}}_{ll} \\
&\overset{(27)}{=} \mathcal{S}_* \widehat{\underline{r}}_{ll} + R\{\mathcal{P}_{\mathrm{add}}(\mathcal{S}_* + R\mathcal{I})^{-1}\mathcal{S}_* - \mathcal{P}_{\mathrm{add}}(\mathcal{S}_* + R\mathcal{I})^{-1}\mathcal{S}_*\}\widehat{\underline{r}}_{ll} = \mathcal{S}_* \widehat{\underline{r}}_{ll}.
\end{aligned}
$$

While $\widehat{\underline{r}}_{ll}$ may be ambiguous, $\mathcal{S}_* \widehat{\underline{r}}_{ll}$ and (if $\Lambda_{*,R}^{-1}$ exists) $\mathcal{P}_{*,R}\widehat{\underline{r}}_{ll}$ are unique.

PROOF OF PROPOSITION 3. In (14), $\mathcal{P}_{*,R}\widehat{\underline{r}}_{ll} \in \mathcal{F}_{\mathrm{add}}$. In order to show that the other term is orthogonal to $\mathcal{F}_{\mathrm{add}}$, we prove that

$$
\mathcal{P}_{\mathrm{add}}(\mathcal{S}_* + R\mathcal{I})^{-1}\mathcal{S}_*(\mathcal{I} - \mathcal{P}_{*,R}) = 0.
$$

This holds because $(\mathcal{S}_* + R\mathcal{I})^{-1}\mathcal{S}_* = \mathcal{I} - (\mathcal{S}_* + R\mathcal{I})^{-1}R$ and (27). $\quad\square$

PROOF OF PROPOSITION 4. The orthogonal projection from $\mathcal{F}$ to $\mathcal{F}_{\mathrm{add}}$ is the same for $\|\cdot\|_*$ and for $\|\cdot\|_R$:

$$
\underset{\underline{r}_{\mathrm{add}} \in \mathcal{F}_{\mathrm{add}}}{\arg\min} \|\underline{r} - \underline{r}_{\mathrm{add}}\|_* = \underset{\underline{r}_{\mathrm{add}} \in \mathcal{F}_{\mathrm{add}}}{\arg\min} \|\underline{r} - \underline{r}_{\mathrm{add}}\|_R \qquad \forall \, \underline{r} \in \mathcal{F},
$$

because the penalty $R\|(\mathcal{I} - \mathcal{P}_{\mathrm{add}})\mathcal{P}_0(\underline{r} - \underline{r}_{\mathrm{add}})\|_2^2$ does not depend on $\underline{r}_{\mathrm{add}}$. Consequently, we may exchange the two norms when projecting to $\mathcal{F}_{\mathrm{add}}$, and we may simplify nested projections:

$$
\begin{aligned}
\mathcal{P}_* \widehat{\underline{r}}_R &\equiv \underset{\underline{r}_{\mathrm{add}} \in \mathcal{F}_{\mathrm{add}}}{\arg\min} \left\| \underline{r}_{\mathrm{add}} - \underset{\underline{r} \in \mathcal{F}_{\mathrm{full}}}{\arg\min} \|\underline{r} - \underline{r}_Y\|_R \right\|_* \\
&= \underset{\underline{r}_{\mathrm{add}} \in \mathcal{F}_{\mathrm{add}}}{\arg\min} \left\| \underline{r}_{\mathrm{add}} - \underset{\underline{r} \in \mathcal{F}_{\mathrm{full}}}{\arg\min} \|\underline{r} - \underline{r}_Y\|_R \right\|_R \\
&= \underset{\underline{r}_{\mathrm{add}} \in \mathcal{F}_{\mathrm{add}}}{\arg\min} \|\underline{r}_{\mathrm{add}} - \underline{r}_Y\|_R = \underset{\underline{r}_{\mathrm{add}} \in \mathcal{F}_{\mathrm{add}}}{\arg\min} \|\underline{r}_{\mathrm{add}} - \underline{r}_Y\|_* \equiv \widehat{\underline{r}}_{\mathrm{add}}. \quad\square
\end{aligned}
$$

A.0.9. *Proofs for Section 4.3.*

PROOF OF PROPOSITION 5. The proof follows from Gao (2003). In particular, continuity on $\mathbb{R}^d$ of the design density $f$ does not hold. However, Condition MLN:B2$'$ states that $f$ is bounded on $[0,1]^d$. For an *upper bound* of $S_{0,0}(\underline{x})$, we choose some smooth density $\widetilde{f}$ and a constant $c$ with $f(\underline{x}) \leq c\widetilde{f}(\underline{x})$ and add $\lfloor (c-1)n \rfloor$ virtual observations with distribution $(c\widetilde{f} - f)/(c-1)$. The density estimator based on all $\lfloor cn \rfloor$ observations is bounded in probability and $S_{0,0}/c$ is smaller. The boundary adjustments are handled analogously.  □

PROOF OF LEMMA 1.   Because of Lemma 2 and Proposition 5, $\mathcal{S}_*$ and $\mathcal{S}_{\mathrm{add}}^{-1}$ are bounded. In Appendix A.0.8, the claim is already shown for small $R$ and (26). Using (17) in Section 4.4, the proof is straightforward for large $R$. □

A.0.10.  *Continuity in $R = \infty$ and proofs for Section 4.4.*   Continuity in $R = 0$ requires no proof.

PROOF OF LEMMA 2.

*Overview.*   MLN showed in Theorem 1′ that the estimator $\widehat{\underline{r}}_{\mathrm{add}}$ is *unique* with probability tending to 1 for $n \to \infty$. Uniqueness is equivalent to $\|\underline{r}\|_* > 0$ for all $\underline{r} \in \mathcal{F}_{\mathrm{add}}$ with $\|\underline{r}\|_2 = 1$. Using their technique of proof, we may even show that the above norm is *bounded away from zero* with probability tending to 1. In this case $\mathcal{S}_{\mathrm{add}}$ has a continuous inverse with respect to $\|\cdot\|_2$.

*Sketch of proof for Theorem* MLN:1′.   The normal equations for $\widehat{\underline{r}}_{\mathrm{add}}$ are

$$\mathcal{S}_{\mathrm{add}} \widehat{\underline{r}}_{\mathrm{add}} = \mathcal{P}_{\mathrm{add}} \underline{r}_L.$$

Define the matrix $\widehat{\mathbf{M}}_k(\underline{\mathbf{x}})$ which depends only on the one-dimensional data $(Y_i, X_{i,k})$, $i = 1, \ldots, n$:

$$\widehat{\mathbf{M}}_k(\underline{\mathbf{x}}) = \frac{1}{n} \sum_{i=1}^n \int K_h(\underline{\mathbf{X}}_i, \underline{\mathbf{x}}) \, d\underline{\mathbf{x}}_{-k} Y_i \begin{bmatrix} 1 & \dfrac{X_{i,k} - x_k}{h_k} \\ \dfrac{X_{i,k} - x_k}{h_k} & \left(\dfrac{X_{i,k} - x_k}{h_k}\right)^2 \end{bmatrix}.$$

With probability tending to 1, $\widehat{\mathbf{M}}_k^{-1}(\underline{\mathbf{x}})$ is continuous, and in this case some $\widehat{\tau}$ is obtained by a continuous mapping of $\mathcal{P}_{\mathrm{add}} \underline{r}_L$:

$$\widehat{\underline{r}}_{\mathrm{add}} = \widehat{\mathcal{T}} \widehat{\underline{r}}_{\mathrm{add}} + \widehat{\tau},$$

where $\widehat{\mathcal{T}}$ is some *contraction*. Therefore the solution is unique and $\mathcal{P}_{\mathrm{add}} \underline{r}_L \mapsto \widehat{\underline{r}}_{\mathrm{add}}$ is continuous. Both $\widehat{\mathcal{T}}$ and $\widehat{\tau}$ depend on $\widehat{\mathbf{M}}_k^{-1}$ and the two-dimensional empirical marginal distribution of $\underline{\mathbf{X}}_i$.

Because $\mathcal{P}_{\mathrm{add}} \underline{r}_L$, even when choosing arbitrary values for $\underline{\mathbf{Y}}$, does not occupy $\mathcal{F}_{\mathrm{add}}$, we cannot (yet) conclude that $\mathcal{S}_{\mathrm{add}}^{-1}$ exists and is continuous.

*Definition of $\mathcal{P}_*$.* Let us now examine the orthogonal (with respect to $\|\cdot\|_*$) projection $\mathcal{P}_*$ from $\mathcal{F}_{\text{full}}$ to $\mathcal{F}_{\text{add}}$,

$$\mathcal{P}_* \breve{\underline{r}} = \operatorname*{arg\,min}_{\breve{\underline{r}}_{\text{add}} \in \mathcal{F}_{\text{add}}} \|\breve{\underline{r}} - \breve{\underline{r}}_{\text{add}}\|_*^2.$$

The normal equations are

$$\mathcal{S}_{\text{add}} \breve{\underline{r}}_{\text{add}} = \mathcal{P}_{\text{add}} \mathcal{S}_* \breve{\underline{r}}.$$

By choosing $\hat{\tau}$ appropriately, one can prove that $\mathcal{P}_{\text{add}} \mathcal{S}_* \breve{\underline{r}} \mapsto \breve{\underline{r}}_{\text{add}}$ is continuous. Because of the uniqueness of $\mathcal{P}_*$, the image of $\mathcal{F}_{\text{full}}$ under the mapping $\mathcal{P}_{\text{add}} \mathcal{S}_*$ is equal to $\mathcal{F}_{\text{add}}$ and therefore $\mathcal{S}_{\text{add}}^{-1}$ exists and is continuous.

*Convergence of $\mathcal{S}_{\text{add}}$.* The operator $\mathcal{S}_{\text{add}}$ depends on *bivariate* terms only. Under Conditions MLN:B1, MLN:B2′ and C1, these terms converge to their theoretical counterparts, which depend on the design density $f$. Furthermore, MLN argued that $\widehat{\mathcal{T}}$ is a contraction (with probability tending to 1) because it converges to $\mathcal{T}$, which is a contraction [MLN:(69)].

*Bandwidths larger than $n^{-1/5}$.* The above calculation assumes that $h$ is proportional to $n^{-1/5}$. In the proof of MLN, one piece was the convergence of $\widehat{\mathcal{T}}$ to $\mathcal{T}$, which depends only on the theoretical design density $f$. In case of oversmoothing, variability is reduced and the expected part is not critical, as Condition MLN:B2′ holds also for smoothed $f$. Hence, Condition C1 may be replaced by Condition C1+. $\quad\square$

PROOF OF LEMMA 3. $\|\mathcal{P}_{\text{add}} \underline{r}_L\|_2$ is essentially univariate and therefore $O_P(1)$.

Define $\mathcal{P}_{\text{add}+} : \mathcal{F}_{\text{full}} \to \mathcal{F}_{\text{add}}$ via $(\mathcal{P}_{\text{add}+} \underline{r})^0(\underline{\mathbf{x}}) = \sum_{k=1}^{d} \int r^0(\underline{\mathbf{x}}) \, d\underline{\mathbf{x}}_{-k}$ and $(\mathcal{P}_{\text{add}+} \underline{r})^k(\underline{\mathbf{x}}) = \int r^k(\underline{\mathbf{x}}) \, d\underline{\mathbf{x}}_{-k}$. By construction, $\mathcal{P}_{\text{add}+}$ is *monotone*: if $r^\ell(\underline{\mathbf{x}}) \geq \breve{r}^\ell(\underline{\mathbf{x}})$ ($\forall \ell, \underline{\mathbf{x}}$), then $(\mathcal{P}_{\text{add}+} \underline{r})^\ell(\underline{\mathbf{x}}) \geq (\mathcal{P}_{\text{add}+} \breve{\underline{r}})^\ell(\underline{\mathbf{x}})$. Note that $\mathcal{P}_{\text{add}}$ is not monotone. Denote by $\mathcal{D}_R$ the contraction $\mathcal{I} - (\mathcal{I} + R^{-1} \mathcal{S}_*)^{-1}$ which is a *pointwise* linear transformation with $\|\mathcal{D}_R\|_{2,\sup} \leq R^{-1} \|\mathcal{S}_*\|_{2,\sup}$.

We want to prove that $\mathcal{P}_{\text{add}+} \mathcal{D}_R \underline{r}_L$ is arbitrarily small when $R^{-1} \mathcal{S}_*$ is small enough. Let us sketch the proof in a simplified case. If $\mathcal{D}_R$ were *diagonal*, we would use the *monotonicity* of $\mathcal{P}_{\text{add}+}$ to obtain an upper bound by replacing $r_L^\ell(\underline{\mathbf{x}})$ by its absolute value ($\forall \ell, \underline{\mathbf{x}}$) and enlarging $\mathcal{D}_R$ to $\|\mathcal{D}_R\|_{2,\sup} \mathcal{I}$:

(28) $$\|(\mathcal{P}_{\text{add}+} \mathcal{D}_R \underline{r}_L)(\underline{\mathbf{x}})\|_2 \leq \|\mathcal{D}_R\|_{2,\sup} \|(\mathcal{P}_{\text{add}+} \underline{r}_{|L|})(\underline{\mathbf{x}})\|_2,$$

where $\underline{r}_{|L|}$ denotes $\underline{r}_L$ with absolute values. The *pointwise* upper bound for $\|(\mathcal{P}_{\text{add}+} \underline{r}_L)(\underline{\mathbf{x}})\|_2$ remains valid if $Y_i$ and $X_{i,k} - x_k$ are replaced by $|Y_i|$ and $|X_{i,k} - x_k|$, respectively.

In practice, positivity of all components is generally not preserved under multiplication by $\mathcal{D}_R$. By Condition MLN:B1, the kernel $K$ has compact support $[-C_1, C_1]$ and therefore the slope terms of $\underline{r}_L$ are bounded by the intercept, $r_L^k(\underline{\mathbf{x}}) \le C_1 r_L^0(\underline{\mathbf{x}})$. The norm $\|(\mathcal{D}_R \underline{r}_L)(\underline{\mathbf{x}})\|_2$ is bounded pointwise in $\underline{\mathbf{x}}$ by $\sqrt{1 + C_1^2 d} \|\mathcal{D}_R\|_{2,\sup} \frac{1}{n} \sum_{i=1}^n K_h(\underline{\mathbf{X}}_i, \underline{\mathbf{x}}) |Y_i|$.

Let $\underline{r}_{|L|}(\underline{\mathbf{x}}) = \sqrt{1 + C_1^2 d} \frac{1}{n} \sum_{i=1}^n K_h(\underline{\mathbf{X}}_i, \underline{\mathbf{x}}) |Y_i| (1, \ldots, 1)$ and (28) holds. Again, $\mathcal{P}_{\text{add}+\underline{r}_{|L|}}(\underline{\mathbf{x}})$ depends on univariate terms only, which is essential in high dimensions, where $\frac{1}{n} \sum_{i=1}^n K_h(\underline{\mathbf{X}}_i, \underline{\mathbf{x}})$ is not of constant order.

We conclude that $\sup_{\underline{\mathbf{x}}} \|(\mathcal{P}_{\text{add}+} \mathcal{D}_R \underline{r}_L)(\underline{\mathbf{x}})\|_2 = O_P(R^{-1} \|\mathcal{S}_*\|_{2,\sup})$.  $\square$

PROOF OF LEMMA 4.  Let $a = \underline{r}_Y - \widehat{\underline{r}}_R$, $b = \widehat{\underline{r}}_R - \widehat{\underline{r}}_{\text{add}}$ and $c = a + b = \underline{r}_Y - \widehat{\underline{r}}_{\text{add}}$. Because $\widehat{\underline{r}}_R$ minimizes $\|\underline{r} - \underline{r}_Y\|_R^2$ and $\|\widehat{\underline{r}}_{\text{add}} - \underline{r}_Y\|_R$ is independent of $R$,

$$(29) \qquad\qquad \|a\|_R \le \|c\|_R = \|c\|_*.$$

Obtain a bound for $\|(\mathcal{I} - \mathcal{P}_{\text{add}}) \widehat{\underline{r}}_R\|_2^2$ using (4):

$$R\|(\mathcal{I} - \mathcal{P}_{\text{add}}) \widehat{\underline{r}}_R\|_2^2 \overset{(4)}{=} \|a\|_R^2 - \|a\|_*^2 \overset{(29)}{\le} \|c\|_*^2 - \|a\|_*^2 = \langle c - a, c + a \rangle_*$$
$$= \langle b, 2c - b \rangle_* \le 2\langle b, c \rangle_* \le 2\|b\|_* \|c\|_*.$$

Then

$$\|(\mathcal{I} - \mathcal{P}_{\text{add}}) \widehat{\underline{r}}_R\|_2^2 \le 2R^{-1}(\sqrt{\|\mathcal{S}_*\|_{2,\sup}} \|\widehat{\underline{r}}_R - \widehat{\underline{r}}_{\text{add}}\|_2) \|\underline{r}_Y - \widehat{\underline{r}}_{\text{add}}\|_*.$$

Because $S_{0,0}(\underline{\mathbf{x}})$ is a density estimate, $\|\mathcal{S}_*\|_{2,\sup} \ge 1$ and we omit the square root.

Because $\widehat{\underline{r}}_{\text{add}}$ is a minimizer, $\|\underline{r}_Y - \widehat{\underline{r}}_{\text{add}}\|_*^2 \le \|\underline{r}_Y\|_*^2 = \frac{1}{n} \sum_{i=1}^n Y_i^2$. For increasing $n$, this is $O_P(1)$ because of Condition MLN:B3$'$.  $\square$

PROOF OF THEOREM 2.  As the expected part is additive [up to $o_P(h^2)$ terms], the nonadditive part consists only of the variance terms: $\|(\mathcal{I} - \mathcal{P}_{*,R}) \widehat{\underline{r}}_{ll}\|_2^2 = O_P(\frac{1}{nh^d})$. By Proposition 4, $\widehat{\underline{r}}_{\text{add}} - \mathcal{P}_{\text{add}} \widehat{\underline{r}}_R = \mathcal{P}_*(\mathcal{I} - \mathcal{P}_{\text{add}}) \widehat{\underline{r}}_R$, where $\mathcal{P}_*$ is continuous (with probability $\to 1$). By Proposition 3, $(\mathcal{I} - \mathcal{P}_{\text{add}}) \widehat{\underline{r}}_R = (\mathcal{S}_* + R\mathcal{I})^{-1} \mathcal{S}_*(\mathcal{I} - \mathcal{P}_{*,R}) \widehat{\underline{r}}_{ll}$ and the claim follows from $\|(\mathcal{S}_* + R\mathcal{I})^{-1} \mathcal{S}_*\|_{2,\sup} \le R^{-1} \|\mathcal{S}_*\|_{2,\sup}$ and uniform continuity of $\mathcal{P}_{*,R}$ (Lemma 1).  $\square$

A.0.11. *Model selection by AIC* (*Section* 4.5).

PROOF OF LEMMA 5.  We use a formula from Rao and Kleffe [(1988), pages 31ff],

$$(30) \qquad \text{cov}(\underline{\varepsilon}^\top \mathbf{B} \underline{\varepsilon}, \underline{\varepsilon}^\top \mathbf{C} \underline{\varepsilon}) = 2\sigma^4 \text{tr}(\mathbf{B}\mathbf{C}) + \kappa \sigma^4 \text{tr}(\mathbf{B} \text{diag}(\mathbf{C})),$$

where $\mathbf{B}$, $\mathbf{C}$ are *symmetric* matrices of dimension $n$ and $\mathbb{E}[\varepsilon_i^4] = (3 + \kappa)\sigma^4$. Using $\mathbf{B} = \mathbf{C} = \frac{1}{2}(\mathbf{M}_R^\top + \mathbf{M}_R)$ we obtain

$$\text{var}(\langle \underline{\varepsilon}, \mathbf{M}_R \underline{\varepsilon} \rangle) = \sigma^4\big(\text{tr}(\mathbf{M}_R \mathbf{M}_R) + \text{tr}(\mathbf{M}_R^\top \mathbf{M}_R) + \kappa \, \text{tr}\,(\text{diag}(\mathbf{M}_R)^2)\big).$$

Note that $\|\mathbf{M}_R\|_{HS}^2 = \text{tr}(\mathbf{M}_R^\top \mathbf{M}_R)$ is known as a Hilbert–Schmidt norm and $\text{tr}(\mathbf{M}_R \mathbf{M}_R) = \langle \mathbf{M}_R^\top, \mathbf{M}_R \rangle_{HS}$ is bounded by $\|\mathbf{M}_R\|_{HS} \|\mathbf{M}_R^\top\|_{HS}$. Hence $\text{var}(\langle \underline{\varepsilon}, \mathbf{M}_R \underline{\varepsilon} \rangle) \leq (2 + |\kappa|)\sigma^4 \, \text{tr}(\mathbf{M}_R^\top \mathbf{M}_R)$ and $\mathbb{E}[\|\mathbf{M}_R \underline{\varepsilon}\|^2] = \sigma^2 \, \text{tr}(\mathbf{M}_R^\top \mathbf{M}_R)$.

Analogously, using $\text{tr}((\mathbf{M}_R^\top \mathbf{M}_R)^2) \leq \|\mathbf{M}_R^\top \mathbf{M}_R\|_{\sup} \text{tr}(\mathbf{M}_R^\top \mathbf{M}_R)$, the variance of $\|\mathbf{M}_R \underline{\varepsilon}\|^2$ is smaller than $(2 + |\kappa|)\sigma^2 \|\mathbf{M}_R^\top \mathbf{M}_R\|_{\sup} \mathbb{E}[\|\mathbf{M}_R \underline{\varepsilon}\|^2]$.  □

A.0.12. *Proofs for Section* 4.6.

PROOF OF PROPOSITION 7. It is well known that $\mathbb{E}[\|\mathbf{M}_{ll} \underline{\varepsilon}\|^2] \propto \frac{1}{nh^d}$ and $\mathbb{E}[\|\mathbf{M}_{\text{add}} \underline{\varepsilon}\|^2] \propto \frac{d}{nh}$. Let us start with the $\lambda^2$ terms: Because $\frac{1}{n\sigma^2}\|\mathbf{M}_{\text{add}} \underline{\varepsilon}\|^2$ is of smaller order than $\frac{1}{n\sigma^2}\|\mathbf{M}_{ll} \underline{\varepsilon}\|^2$, the mixed term $\frac{2}{n\sigma^2}\langle \mathbf{M}_{ll} \underline{\varepsilon}, \mathbf{M}_{\text{add}} \underline{\varepsilon} \rangle$ is bounded by the Cauchy–Schwarz inequality. Furthermore, $(1 - \mathbb{E})[\frac{1}{n\sigma^2}\|\mathbf{M}_{ll} \underline{\varepsilon}\|^2]$ is negligible compared to $\mathbb{E}[\frac{1}{n\sigma^2}\|\mathbf{M}_{ll} \underline{\varepsilon}\|^2]$ (see Lemma 5), indicating that the *inverse* of the $\lambda^2$ terms is $O_P(nh^d)$.

For the $\lambda$-linear terms, it is obvious that $\frac{2}{n\sigma^2}(1 - \mathbb{E})[\langle \underline{\varepsilon}, \mathbf{M}_{ll} \underline{\varepsilon} \rangle]$ is the largest stochastic term. It remains to show that $\frac{2}{n\sigma^2}\mathbb{E}[\langle \mathbf{M}_{ll} \underline{\varepsilon}, \mathbf{M}_{\text{add}} \underline{\varepsilon} \rangle]$ is nonnegative, as a nonnegative coefficient of $\lambda$ indicates that the minimum is at $\lambda_{\min} \leq 0$ ($R = \infty$). As we are using a fixed uniform design, the local linear and the Nadaraya–Watson estimator coincide (ignoring the boundary). Hence, we are interested in the covariance of a multivariate and a univariate Nadaraya–Watson estimator with *nonnegative kernel weights*, whose hat matrices have therefore nonnegative elements.  □

**Acknowledgments.** The authors would like to thank the editors and referees for their helpful comments that led to a considerable improvement of this paper.

DEPARTMENT OF BIOSTATISTICS, ISPM
UNIVERSITY OF ZURICH
SUMATRASTRASSE 30
CH-8006 ZURICH
SWITZERLAND
E-MAIL: mstuder@ifspm.unizh.ch